\documentclass[11pt]{article}

\usepackage{graphicx}
\usepackage{nameref}
\usepackage[shortlabels,inline]{enumitem}
\usepackage{empheq}
\usepackage[shortlabels]{enumitem}
\setlist[enumerate]{nosep}
\usepackage[doc]{optional}
\usepackage{xcolor}
\usepackage[colorlinks=true,
            linkcolor=refkey,
            urlcolor=lblue,
            citecolor=red]{hyperref}
\usepackage{float}
\usepackage{soul}
\usepackage{graphicx}
\definecolor{labelkey}{rgb}{0,0.08,0.45}
\definecolor{refkey}{rgb}{0,0.6,0.0}
\definecolor{Brown}{rgb}{0.45,0.0,0.05}
\definecolor{lime}{rgb}{0.00,0.8,0.0}
\definecolor{lblue}{rgb}{0.5,0.5,0.99}

 \usepackage{mathpazo}

\colorlet{hlcyan}{cyan!30}

\usepackage{stmaryrd}
\usepackage{amssymb}
\makeatletter
\def\namedlabel#1#2{\begingroup
   \def\@currentlabel{#2}%
   \label{#1}\endgroup
}
\makeatother

\newcommand{\sepp}{\setlength{\itemsep}{-2pt}}

\newcommand{\bb}{\ensuremath{\mathbf{b}}}
\newcommand{\bA}{\ensuremath{\mathbf{A}}}
\newcommand{\bt}{\ensuremath{\mathbf{t}}}

\newcommand{\bx}{\ensuremath{\mathbf{x}}}

\newcommand{\menge}[2]{\big\{{#1}~\big |~{#2}\big\}}

\newcommand{\To}{\ensuremath{\rightrightarrows}}

\newcommand{\fenv}[1]%
{\ensuremath{\,\overrightarrow{\operatorname{env}}_{#1}}}
\newcommand{\benv}[1]%
{\ensuremath{\,\overleftarrow{\operatorname{env}}_{#1}}}

\newcommand{\scal}[2]{\left\langle{#1},{#2}  \right\rangle}

\newcommand{\RR}{\ensuremath{\mathbb R}}

\newcommand{\RPP}{\ensuremath{\mathbb{R}_{++}}}

\newcommand{\RX}{\ensuremath{\,\left]-\infty,+\infty\right]}}
\newcommand{\RRX}{\ensuremath{\,\left[-\infty,+\infty\right]}}

\newcommand{\ZZ}{\ensuremath{\mathbf Z}}
\newcommand{\NN}{\ensuremath{\mathbb N}}

\newcommand{\dom}{\ensuremath{\operatorname{dom}}}
\newcommand{\argmin}{\ensuremath{\operatorname{argmin}}}
\newcommand{\prox}{\ensuremath{\operatorname{Prox}}}
\newcommand{\reli}{\ensuremath{\operatorname{ri}}}
\newcommand{\inte}{\ensuremath{\operatorname{int}}}

\newcommand{\closu}{\ensuremath{\operatorname{cl}}}

\newcommand{\ran}{\ensuremath{{\operatorname{ran}}\,}}
\newcommand{\zer}{\ensuremath{\operatorname{zer}}}

\newcommand{\Fix}{\ensuremath{\operatorname{Fix}}}
\newcommand{\Id}{\ensuremath{\operatorname{Id}}}

\newcommand{\bX}{\ensuremath{{\mathbf{X}}}}

\newcommand{\bC}{\ensuremath{{\mathbf{C}}}}

\newcommand{\bR}{\ensuremath{{\mathbf{R}}}}

\newcommand{\bG}{\ensuremath{{\mathbf{G}}}}

\newcommand{\by}{\ensuremath{\mathbf{y}}}
\newcommand{\bz}{\ensuremath{\mathbf{z}}}
\newcommand{\bg}{\ensuremath{\mathbf{g}}}

\newcommand{\dd}{\ensuremath{\mathbf{d}}}

\newcommand{\dper}{\ensuremath{\mathbf{\Delta}^{\perp}}}
\newcommand{\ff}{\ensuremath{\mathbf{f}}}

\newcommand{\GF}{\ensuremath{\Gamma_{0}(X)}}

{\begin{list}{}{%
\settowidth{\labelwidth}{\textrm{#1~}}%
\setlength{\leftmargin}{\labelwidth+\labelsep}}}
{\end{list}}
\usepackage{amsthm}
\usepackage[capitalize,nameinlink]{cleveref}
\crefname{equation}{}{equations}
\crefname{chapter}{Appendix}{chapters}
\crefname{item}{}{items}
\crefname{enumi}{}{}
\newtheorem{theorem}{Theorem}[section]
\newtheorem{lemma}[theorem]{Lemma}

\newtheorem{corollary}[theorem]{Corollary}

\newtheorem{proposition}[theorem]{Proposition}

\newtheorem{example}[theorem]{Example}
\newtheorem{ex}[theorem]{Example}
\newtheorem{fact}[theorem]{Fact}
\newtheorem{remark}[theorem]{Remark}

\newcommand{\bv}{\ensuremath{{\mathbf{v}}}}

\providecommand{\RR}{\mathbb{R}}

\providecommand{\ran}{\operatorname{ran}}

\providecommand{\dom}{\operatorname{dom}}

\providecommand{\epi}{\operatorname{epi}}

\providecommand{\Id}{\operatorname{{ Id}}}

\providecommand{\ZZ}{{\bf Z}}

\providecommand{\dD}{\ensuremath{{\mathbf{\Delta}}}}
\providecommand{\bu}{{\bf u}}
\providecommand{\bv}{{\bf v}}

\providecommand{\argmin}{\mathrm{arg}\!\min}

\providecommand{\To}{\rightrightarrows}

\providecommand{\NN}{\mathbb{N}}

\providecommand{\ran}{\operatorname{ran}}

\providecommand{\Id}{\operatorname{Id}}

\providecommand{\zer}{\operatorname{zer}}

\providecommand{\RR}{\mathbb{R}}
\providecommand{\NN}{\mathbb{N}}

\providecommand{\bt}{b}

\definecolor{myblue}{rgb}{.8, .8, 1}

\allowdisplaybreaks 

\begin{document}

\title{\textsc{
Attouch-Th\'era Duality, Generalized Cycles and Gap Vectors
}}

\author{
Salihah Alwadani\thanks{
Mathematics, University
of British Columbia,
Kelowna, B.C.\ V1V~1V7, Canada. E-mail:
\texttt{saliha01@mail.ubc.ca}.},~
Heinz H.\ Bauschke\thanks{
Mathematics, University
of British Columbia,
Kelowna, B.C.\ V1V~1V7, Canada. E-mail:
\texttt{heinz.bauschke@ubc.ca}.},~
and
Xianfu Wang\thanks{
Mathematics, University
of British Columbia,
Kelowna, B.C.\ V1V~1V7, Canada. E-mail:
\texttt{shawn.wang@ubc.ca}.}
}

\date{January 13, 2021 (Revision April 22, 2021)}
\maketitle

\vskip 8mm

\begin{abstract} \noindent Using the Attouch-Th\'era duality, we study the cycles, gap vectors and
fixed point sets of compositions of proximal mappings. Sufficient conditions are given for the existence of
cycles and gap vectors. A primal-dual framework provides an exact relationship between
the cycles and gap vectors. We also introduce the generalized cycle and gap vectors to tackle the case
when the classical ones do not exist. Examples are given to illustrate our results.
\end{abstract}

{

}

\section{Introduction}


Throughout this paper, we assume that $X$ is a real Hilbert space with inner product
    $\scal{\cdot}{\cdot}\colon X\times X\to\RR$
and induced norm $\|\cdot\|=\sqrt{\scal{\cdot}{\cdot}}$. The set of proper lower semicontinuous
convex functions from $X$ to $]-\infty,+\infty]$  is denoted by
$\Gamma_{0}(X)$. In the product space $\bX=X^m$ with $m\in\NN=\{1, 2, \dots\}$, we let
$$ \Delta=\menge{(x,\ldots,x)}{x\in X}, $$
$$\bR:\bX\rightarrow\bX:(x_{1},x_{2},\ldots, x_{m})\mapsto (x_{m},x_{1},\ldots, x_{m-1}), \text{ and }$$
$$\|\bx\|=\sqrt{\scal{\bx}{\bx}}=\sqrt{\|x_{1}\|^2+\cdots+\|x_{m}\|^2}$$
where $\bx=(x_{1},x_{2},\ldots, x_{m}).$
For a finite family of functions $(f_{i})_{i=1}^{m}$ in $\GF$, define its separable sum by
\begin{equation}\label{e:dsum}
\ff=f_{1}\oplus\cdots\oplus f_{m}: \bX\rightarrow\RX: (x_{1},\ldots, x_{m})\mapsto \sum_{i=1}^{m}f_{i}(x_{i}).
\end{equation}
The proximal mapping  of $f_{i}$ is defined by
$\prox_{f_{i}}=(\Id+\partial f_{i})^{-1}$ where $\partial f_{i}$ denotes the subdifferential of $f_{i}$.
A \emph{classical cycle} of $\ff$ is a vector $\bz=(z_{1},\ldots,z_{m})\in \bX$ such that
\begin{subequations}
\begin{equation}\label{e:cycle}
z_{1} =\prox_{f_{1}}z_{m},  \quad z_{2}=\prox_{f_{2}}z_{1}, \quad z_{3}=\prox_{f_{3}}z_{2}, \cdots,
\end{equation}
\begin{equation}\label{e:cycle:end}
z_{m-1} =\prox_{f_{m-1}}z_{m-2},  \quad z_{m}=\prox_{f_{m}}z_{m-1}.
\end{equation}
\end{subequations}
The set of all classical cycles of $\ff$ will be denoted by $\ZZ$.
Since $\partial\ff=\partial f_{1}\times\cdots\times \partial f_{m}$ and
$\prox_{\ff}=(\prox_{f_{1}},\cdots,\prox_{f_{m}})$,
in the framework of product space $\bX$, with $\bz=(z_{1},\ldots, z_{m})$,
the operator form of
\eqref{e:cycle}--\eqref{e:cycle:end} is
\begin{equation}\label{e:ab:cycle}
\bz=\prox_{\ff}\bR \bz, \text{ equivalently, }
\end{equation}
in terms of monotone operators
\begin{equation}\label{e:sub:form}
0\in \partial \ff(\bz)+\bz-\bR\bz.
\end{equation}
Here the displacement mapping $\Id-\bR$ is maximally monotone but not a gradient of convex function
unless $m=2$.
See \cite[Section VI-B]{comb20} and
\cite[Section 4]{com-pes} for proximal cycles from the perspective of Nash equilibrium.

Many authors have studied cycles and asymptotic behaviors of compositions of proximal mappings or
resolvents (more generally, averaged mappings); see, e.g.,
\cite{abjw2020,bbr78, baillon, baillon1,  BC2017, Victoria, reich, comb20, com-pes, GR, reich1, reich2, wangb}.
For the compositions of two proximal mappings or resolvents, the investigation on the geometry of
fixed point sets and gap vectors has matured; see
\cite{BCR,wangb}. For the compositions of three or more proximal mappings or resolvents,
when the intersection of fixed point sets of proximal mapping or resolvents is empty,
the geometry of fixed point sets and gap vectors still await to be explored.
On the one hand, a systematic study of cycles and gap vectors for compositions of
proximal mappings does not exist in the literature;
on the other hand,
it is not clear what one should do when the cycles and gap vectors do not exist.
In \cite{abjw2020}, we carried out the study of extended gap vectors for projection mappings and provided an answer to
the geometry conjecture, which concerns the situation where
the classical gap vector does not necessarily exist.

\emph{The goal of this paper is to give a systematic study of the cycles, gap vectors and
fixed points of compositions of proximal mappings, and to provide generalized cycles and
gap vectors when the classical ones do not necessarily
exist. Our investigation relies to a large extent on the
Attouch-Th\'era duality and convex analysis.}

The remainder of the paper is organized as follows. Section~\ref{s:preliminary}
collects some preliminary results and facts on convex functions and monotone operators needed in the sequel.
In Section~\ref{s:cycle:fixed} we provide a systematic study of classical cycles, gap vector, and
fixed point sets of cyclic proximal mappings. It turns out that the fixed point sets of
a cyclic composition
of proximal mappings is just a shift of the other via gap vectors. That is, for
$F_{i}:=\Fix\big(\prox_{f_{i}}\cdots\prox_{f_{1}}\prox_{f_{m}}\cdots\prox_{f_{i+1}}\big)$ with
$F_{m+1}:=F_{1}$,
one has $F_{i+1}=F_{i}-y_{i+1}$ for $i=1,\ldots, m$, where $\by=(y_{1},\ldots, y_{m})$ is the gap vector.
In Section~\ref{s:gcycle:gfixed}, using the Attouch-Th\'era duality
we investigate the generalized cycle and gap vectors, and show that they can recover the classical ones
whenever the latter exist.
Section~\ref{s:lowerh:exact} is devoted to conditions under which $\ff\Box\iota_{\Delta}$ is lower
semicontinuous, in turn, these guarantee that the classical cycles and gap vectors exist.
Section~\ref{s:example} gives two examples to illustrate our results on the generalized cycle and gap vector.
Finally, Section~\ref{s:compute} shows how to use the forward-backward iteration scheme to compute the
generalized cycle and gap vectors.

\noindent {\bf Notation}. We follow the notation of \cite{BC2017}, where one will find a detailed account of the following notions.
Let $f, g:X\rightarrow\RX$.
The Fenchel conjugate of $f$ is
$$f^*:X\rightarrow\RRX:x^*\mapsto\sup_{x\in X}(\scal{x}{x^*}-f(x)).$$ The infimal convolution
of $f, g$  is $ f\Box g:X\rightarrow\RRX:x\mapsto\inf_{y\in X}(f(y)+g(x-y)),$
and it is exact at a point $x\in X$ if
$(\exists y\in X)\ (f\Box g)(x)=f(y)+g(x-y);$
$f\Box g$ is exact if it is exact at every point of its domain.
The subdifferential of $f$ is the set-valued operator
$$\partial f:X\To X: x\mapsto\menge{x^*\in X}{(\forall y\in X)\ f(y)\geq f(x)+\scal{x^*}{y-x}}.$$
We use $\closu f$ for the lower semicontinuous hull of $f$.
For a set $C\subset X$, its indicator function is defined by
$$\iota_{C}(x)=\begin{cases}
0, & \text{ if $x\in C$,}\\
+\infty, & \text{ if $x\not\in C$.}
\end{cases}
$$
The relative interior, interior, and closure of $C$ will be denoted by $\reli C$, $\inte C$, and $\overline{C}$
respectively.
When the set $C$ is nonempty closed convex, we
write $P_{C}=\prox_{\iota_{C}}$ for the projection operator and $N_{C}=\partial \iota_{C}$
for the normal cone.
An operator $N:X\rightarrow X$ is nonexpansive if $(\forall x, y \in X)\ \|Nx-Ny\|\leq \|x-y\|$;
firmly nonexpansive if $2N-\Id$ is nonexpansive; $\beta$-cocercive if $\beta N$ is firmly
nonexpansive for some $\beta\in\RPP$. Prime examples of firmly nonexpansive mappings are
proximal mappings of convex functions.
As usual, $\Fix N=\menge{x\in X}{Nx=x}$ denotes the set of fixed points of $N$.
For a monotone operator $A: X \To X$, the sets
 $\dom A =\menge{x\in X}{Ax\neq\varnothing},
 \ran A =\menge{u\in X}{(\exists x \in X) \ u\in Ax}$, $\ker A=\zer(A)=A^{-1}(0)$
 are the domain, range, and kernel or zeros of $A$ respectively.
 It will be convenient to write
$\widetilde{A} =(-\Id)\circ A^{-1}\circ (-\Id).$
Let $(x_{n})_{n\in\NN}$ be a sequence in $X$. We will write
$x_{n}\rightarrow x$ if $(x_{n})_{n\in\NN}$ converges strongly to $x$, and
$x_{n}\rightharpoonup x$ if $(x_{n})_{n\in\NN}$ converges weakly to $x$.
 For valuable references on convex functions and monotone operators,
see, e.g., \cite{Rock70,Rock98,Zalinescu,Simons2}.

\noindent {\bf Blanket assumptions}. Throughout the paper, we shall assume that
\begin{enumerate}
\item $(f_{i})_{i=1}^{m}$ are in $\GF$, and $\ff$ is given by \eqref{e:dsum}.
\item \begin{equation}\label{e:conj:proper}
\dom (\ff^*+\iota_{\Delta}^*)=\dom (\ff^*+\iota_{\dper})\neq \varnothing, \text{equivalently, $\dom \ff^*\cap \dper\neq\varnothing$.}
\end{equation}
This will assure that $\ff\Box\iota_{\Delta}$
is proper convex, and possess a continuous affine minorant, see, e.g., \cite[Proposition 13.12(ii)]{BC2017}.
\end{enumerate}

\section{Auxiliary results}\label{s:preliminary}
In this section we give some conditions under which $\dom \ff^*\cap \dper\neq\varnothing$, subdifferential
properties concerning the lower semicontinuous hulls and infimal convolutions, and the Attouch-Th\'era
duality for pairs of monotone operators.
\subsection{Properties of $\ff^*+\iota_{\Delta}^*, \ff\Box\iota_{\Delta},$ and $\closu(\ff\Box\iota_{\Delta})$}
Because $\ff$ is separable, we have $\ff^*=f_{1}^*\oplus\cdots \oplus f_{m}^*.$
Recall that a function $g\in\GF$ is \emph{coercive}
if $\lim_{\|x\|\rightarrow\infty}g(x)=+\infty$,
and \emph{supercoercive} if
$\lim_{\|x\|\rightarrow\infty}(g(x)/\|x\|)=+\infty.$

\begin{lemma} Suppose that one of the following hold:
\begin{enumerate}
\item \label{i:affine}
There exists $\bx^*=(x_{1}^*,\ldots, x_{m}^*)\in\dper$ and $\mu_{1},\ldots,\mu_{m}\in\RR$ such that
$(\forall 1\leq i\leq m)(\forall x_{i}\in X)\  f_{i}(x_{i})\geq\scal{x_{i}^*}{x_{i}}-\mu_{i}$.
\item \label{i:one:bound} One of the functions $f_{i}$ is supercoercive, e.g., having a bounded domain.
\item\label{i:bounded} $\inf f_{i}>-\infty$ for every $i=1,\ldots, m$.
\end{enumerate}
Then \eqref{e:conj:proper} holds.
\end{lemma}
\begin{proof}
Because $(\ff\Box\iota_{\Delta})^*=\ff^*+\iota_{\Delta}^*$,
\eqref{e:conj:proper} holds if and only if $\ff\Box\iota_{\Delta}$ possesses a continuous minorant.

Suppose \ref{i:affine} holds.
Then
$$(\forall 1\leq i\leq m)(\forall x_{i}\in X)\ f_{i}(x_{i}-d)\geq \scal{x_{i}^*}{x_{i}-d}-\mu_{i}=\scal{x_{i}^*}{x_{i}}-\scal{x_{i}^*}{d}-\mu_{i}.$$
 Summing up these $m$-inequalities gives
 $$\sum_{i=1}^{m}f_{i}(x_{i}-d)\geq \sum_{i=1}^{m}\scal{x_{i}^*}{x_{i}}-\sum_{i=1}^{m}\mu_{i}=\scal{\bx^*}{\bx}-\sum_{i=1}^{m}\mu_{i}.$$
 Taking infimum over $d\in X$ gives
 $(\ff\Box\iota_{\Delta})(\bx)\geq \scal{\bx^*}{\bx}-\sum_{i=1}^{m}\mu_{i},$
 as required.

Suppose that \ref{i:one:bound} holds, say $f_{m}$ being supercoercive. For each $f_{i}$ with $1\leq i\leq m-1$,
by Bronsted-Rockafellar's theorem \cite[Theorem 16.58]{BC2017}
 we can choose an affine functional $x_{i}^*\in X$, $\mu_{i}\in\RR$ such that
$$(\forall x_{i}\in X)\  f_{i}(x_{i})\geq\scal{x_{i}^*}{x_{i}}-\mu_{i}.$$
Put $x_{m}^*=-(x_{1}^*+\cdots+x_{m-1}^*)$. Since $f_{m}$ is supercoercive, the function
$f_{m}-\scal{x_{m}^*}{\cdot}$ is coercive, so it attains its minimum at $-\mu_{m}\in \RR$. Then
$f_{m}\geq \scal{x_{m}^*}{\cdot}-\mu_{m}$. Hence \ref{i:affine} applies.

 \ref{i:bounded} is a special case of \ref{i:affine} with $x_{i}^*=0$ for $i=1,\ldots, m$.
\end{proof}

Some elementary properties of $\ff\Box\iota_{\Delta}$ and its lower semicontinuous hull
$\closu(\ff\Box\iota_{\Delta})$ come next.
\begin{lemma}\label{l:range}
 We have
$\ran\partial[\closu(\ff\Box\iota_{\Delta})]\subseteq\dper.$
\end{lemma}
\begin{proof} Let $\bx^*\in \partial[\closu(\ff\Box\iota_{\Delta})](\bx)$. Then
$(\forall \bz\in \bX) \ [\closu(\ff\Box\iota_{\Delta})](\bz)\geq [\closu(\ff\Box\iota_{\Delta})](\bx)+\scal{\bx^*}{\bz-\bx}.$
Let $\bz=\bu+\bv$ with $\bu\in\dom\ff$ and $\bv\in\Delta$. We have
$$\ff(\bu)=\ff(\bu)+\iota_{\Delta}(\bv)\geq [\closu(\ff\Box\iota_{\Delta})](\bx)+\scal{\bx^*}{\bu-\bx}+
\scal{x^*}{\bv}.$$
Since $\Delta$ is a subspace and $\ff(\bu)<+\infty$, we obtain
$\scal{\bx^*}{\bv}=0.$
\end{proof}
\begin{lemma}\label{c:subdiff}
 For every $\dd\in\Delta$ and $\bx\in \bX$, we have
\begin{enumerate}
\item\label{i:ftrans}
$(\ff\Box\iota_{\Delta})(\bx +\dd)=(\ff\Box\iota_{\Delta})(\bx),$
and
$[\closu(\ff\Box\iota_{\Delta})](\bx+\dd)=[\closu(\ff\Box\iota_{\Delta})](\bx).$
\item\label{i:strans}
$\partial(\ff\Box\iota_{\Delta})(\bx +\dd)=\partial (\ff\Box\iota_{\Delta})(\bx), \text{ and }$
$\partial[\closu(\ff\Box\iota_{\Delta})](\bx+\dd)=\partial[\closu(\ff\Box\iota_{\Delta})](\bx).$
\end{enumerate}
\end{lemma}
\begin{proof}
\ref{i:ftrans}: Since $\Delta$ is a subspace, we have
$\dd-\Delta=\Delta$ and $-\Delta=\Delta$.
Then
\begin{align}
(\ff\Box\iota_{\Delta})(\bx+\dd) &=\inf_{\bv\in \Delta}(\ff(\bx+\dd-\bv)+\iota_{\Delta}(\bv))\\
& =\inf_{\bv\in \dd-\Delta}\ff(\bx+\bv)
=\inf_{\bv\in \Delta}\ff(\bx+\bv)=\inf_{-\bv\in \Delta}(\ff(\bx-(-\bv))
=(\ff\Box\iota_{\Delta})(\bx).
\end{align}
\ref{i:strans}: Apply \ref{i:ftrans} and the subdifferential definition.
\end{proof}

\subsection{Subdifferentials}

\begin{lemma}\label{l:closed}
 Let $f:X\rightarrow\RX$ be proper and convex, and $x\in X$. Then the following assertions hold:
 \begin{enumerate}
 \item\label{i:nonemptys} If $\partial f(x)\neq\varnothing$, then $f$ is lower semicontinuous at $x$.
 \item\label{i:equal} If $f(x)=(\closu f)(x)$, that is,
$f$ is lower semicontinuous at $x$, then
$\partial f(x)=\partial (\closu f)(x)$.
\item \label{i:inclusion} In general,
 $\partial f\subseteq\partial (\closu f)$.
\end{enumerate}
\end{lemma}
\begin{proof}
\ref{i:nonemptys}: See \cite[Theorem 2.4.1(ii)]{Zalinescu}.

\ref{i:equal}:
Let $x^*\in \partial (\closu f)(x)$. Then
$(\forall y\in X)\ (\closu f)(y)\geq (\closu f)(x) +\scal{x^*}{y-x}.$
Because $f\geq\closu f$ and $f(x)=(\closu f)(x)$, we have
$(\forall y\in X)\ f(y)\geq f(x) +\scal{x^*}{y-x},$
so $x^*\in\partial f(x).$

Let $x^*\in \partial f(x)$. Then
$(\forall y\in X)\ f(y)\geq f(x) +\scal{x^*}{y-x}.$
Taking the lower semicontinuous envelopes with respect to $y$ both sides, we have
$(\forall y\in X)\ (\closu f)(y)\geq f(x) +\scal{x^*}{y-x}.$
Since $f(x)=(\closu f)(x)$, we obtain $x^*\in\partial\closu f(x)$.

\ref{i:inclusion}: Let $y\in X$. If $\partial f(y)=\varnothing$, then $\partial f(y)\subseteq\partial(\closu f)(y)$.
If $\partial f(y)\neq\varnothing$, then
$f$ is lower semicontinuous at $y$ by \ref{i:nonemptys}. Apply \ref{i:equal}.
\end{proof}

The inclusion in Lemma~\ref{l:closed}~\ref{i:inclusion} can be proper, e.g., $\partial\iota_{C}(x)=\varnothing$ while
$\partial\iota_{\overline{C}}(x)\neq\varnothing$ for $x\in\overline{C}\setminus C$ when
$C$ is convex and $\inte C\neq\varnothing$.

\begin{lemma}\label{l:exact}
Let $f,g\in\GF$ and $x, y\in X$. Then the following assertions hold:
\begin{enumerate}
\item \label{i:exact:funct}
If $(f\Box g)(x)=f(y)+g(x-y)$, then
$\partial (f\Box g) (x)=\partial f(y)\cap \partial g(x-y).$

\item \label{i:exact:sub}
If $\partial f(y)\cap \partial g(x-y)\neq \varnothing$, then
$(f\Box g)(x)=f(y)+g(x-y)$ and
$$\partial (f\Box g) (x)=\partial f(y)\cap \partial g(x-y).$$

\item \label{i:g:sub} In general,
$ \partial (f\Box g) (x) \supseteq \partial f(y)\cap \partial g(x-y).$
\end{enumerate}
\end{lemma}

\begin{proof} \ref{i:exact:funct}: See \cite[Propostition 16.61(i)]{BC2017}.

\ref{i:exact:sub}: When $\partial f(y)\cap \partial g(x-y)\neq \varnothing$,
\cite[Proposition 16.61(ii)]{BC2017} gives $(f\Box g)(x)=f(y)+g(x-y)$.
Apply  \ref{i:exact:funct} to obtain
$\partial (f\Box g) (x)=\partial f(y)\cap \partial g(x-y).$

\ref{i:g:sub}: If $\partial f(y)\cap \partial g(x-y)=\varnothing$, it is clear.
If $\partial f(y)\cap \partial g(x-y)\neq\varnothing$, apply
\ref{i:exact:sub}.
\end{proof}

\subsection{Monotone operators}
The key tool we shall use is the following Attouch-Th\'era duality.
\begin{fact}[Attouch-Th\'era duality \cite{AtTh}]\label{dualityat}
Let $A,B:X\To X$ be maximally monotone operators. Let
$S$ be the solution set of the primal problem:
\begin{equation}\label{theraprimal}
\text{ find $x\in X$ such that } 0\in Ax+Bx.
\end{equation}
Let $S^*$ be the solution set of the dual problem associated with the ordered pair $(A,B)$:
\begin{equation}
\label{theradual1}
\text{ find $x^*\in X$ such that } 0\in A^{-1}x^*+\widetilde{B}(x^*).
\end{equation}
Then
\begin{enumerate}
\item $S=\menge{x\in X}{(\exists\ x^*\in S^*)\ x^*\in Ax \text{
and } -x^*\in Bx}$.
\item $S^*=\menge{x^*\in X}{(\exists\ x\in S)\ x\in A^{-1}x^*
\text{ and } -x\in \widetilde{B}(x^*)}.$
\end{enumerate}
\end{fact}

Important properties of the circular right shift operator come as follows.
\begin{fact}\label{f:shift} For the circular right shift operator $\bR$, the following assertions hold:
\begin{enumerate}
\item\label{i:r1} $\Id-\bR$ is maximally monotone.
\item\label{i:r2} $(\Id-\bR)^{-1}=\frac{1}{2}\Id+N_{\dper}+T$
where $T:\bX\rightarrow \bX$ is a skew operator defined by
$$T=\frac{1}{2m}\sum_{k=1}^{m-1}(m-2k)\bR^{k}.$$ In particular,
$\dom (\Id-\bR)^{-1}=\dper$.
\item\label{i:r3} $(\frac{1}{2}\Id+T)^{-1}=\Id-\bR+2P_{\dD}$.
\end{enumerate}
\end{fact}
\begin{proof}
\ref{i:r1}: See \cite[Theorem 6.1]{bwy14} or \cite[Exercise 12.16]{Rock98}.
\ref{i:r2}: See \cite[Section 2.2]{abjw2020} or \cite{ABRW}. \ref{i:r3}: See \cite[Theorem 2.4]{abjw2020}.
\end{proof}

\begin{fact}\label{l:inv:strong}
 Let $N:X\rightarrow X$ be nonexpansive. Then $(\Id-N)^{-1}$ is
$1/2$-strongly monotone; equivalently, $\Id-N$ is $1/2$-cocoercive.
\end{fact}
\begin{proof} See \cite[Proposition 4.11 and Example 22.7]{BC2017}.
\end{proof}

\section{Classical cycles, gap vectors and fixed point sets of cyclic proximal mappings}\label{s:cycle:fixed}
In this section we give a comprehensive study
on the relationship among the classical cycles, gap vectors and fixed point sets for compositions of proximal mappings.

Define
\begin{equation}\label{e:Fi}
F_{i}:=\Fix\big(\prox_{f_{i}}\cdots\prox_{f_{1}}\prox_{f_{m}}\cdots\prox_{f_{i+1}}\big),\text{ and }
\end{equation}
$$Q_{i}:\bX\rightarrow X: \bz\mapsto z_{i}.$$
Some basic properties and relationship among the $F_{i}$'s, and $\ZZ$ are given below.
\begin{lemma}\label{l:inter:rel}
 For every $1\leq i, j\leq m$, the following assertions hold:
\begin{enumerate}
\item\label{i:intersection} $\bigcap_{i=1}^{m}\argmin f_{i}\subseteq \bigcap_{i=1}^{m} F_{i}.$
In particular, if there exists $F_{i}=\varnothing$, then $\bigcap_{i=1}^{m}\argmin f_{i}=\varnothing$.
\item\label{i:convexset}
 $F_{i}$ is closed and convex, and $F_{i}\subseteq\ran(\prox_{f_{i}})$.
\item \label{i:FZ}
$F_{i}\neq \varnothing$ if and only if $F_{j}\neq\varnothing$ if and only if $\ZZ$ is nonempty.
\item\label{i:zconvex} $\ZZ$ is closed and convex, and
$\ZZ\subseteq F_{1}\times\cdots\times F_{m}$.
\item\label{i:proj} The mapping $Q_{i}:\ZZ\rightarrow F_{i}$ is bijective, in particular,
$Q_{i}(\ZZ)=F_{i}.$

\item\label{i:iterative}
For $1\leq i\leq m-1$,
\begin{equation}\label{e:image}
\prox_{f_{i+1}}(F_{i})= F_{i+1}, \text{ and $\prox_{f_{1}}(F_{m})= F_{1}.$}
\end{equation}
Consequently,
$\prox_{\ff}\bR(F_{1}\times\cdots F_{m})= F_{1}\times\cdots\times F_{m}.$
\item\label{i:twos} When $m=2$, we have
$\prox_{f_{1}}(F_{2})=F_{1}, \prox_{f_{2}}(F_{1})=F_{2}.$
\end{enumerate}
\end{lemma}
\begin{proof}
\ref{i:intersection}: This follows from $\argmin f_{i}=\Fix\prox_{f_{i}}$ for
$1\leq i\leq m$ and
$\bigcap_{i=1}^{m}\Fix\prox_{f_{i}}\subseteq F_{j}$ for $1\leq j\leq m$.

\ref{i:convexset}: Since each $\prox_{f_{i}}$ is firmly nonexpansive, the mapping
$\prox_{f_{i}}\cdots\prox_{f_{1}}\prox_{f_{m}}\cdots\prox_{f_{i+1}}$
is nonexpansive. By \cite[Corollary 4.24]{BC2017}, $F_{i}$ is closed and convex.
The relation $F_{i}\subseteq\ran(\prox_{f_{i}})$ is from the definition of $F_{i}$.

\ref{i:FZ}: This is clear from the definition of $F_{i}, F_{j}, \ZZ$.

\ref{i:zconvex}: Since $\ZZ=\Fix(\prox_{\ff}\bR)$ by \eqref{e:ab:cycle} and $\prox_{\ff}\bR$ is nonexpansive,
it suffices to apply \cite[Corollary 4.24]{BC2017}.
The relation
$\ZZ\subseteq F_{1}\times\cdots\times F_{m}$ is from the definition of $\ZZ$ and $F_{i}$'s.

\ref{i:proj}: It suffices to show $Q_{i}$ is injective. Suppose
$\bz=(z_{1},\ldots, z_m), \tilde{\bz}=(\tilde{z}_{1},\ldots, \tilde{z}_{m})\in\ZZ$ and
$Q_{i}(\bz)=Q_{i}(\tilde{\bz})$, i.e.,
$z_i=\tilde{z}_{i}$. Because $\bz, \tilde{\bz}\in \ZZ$, they are cycles, so we have
$z_{i+1}=\prox_{f_{i+1}}z_{i}=\prox_{f_{i+1}}\tilde{z}_{i}=\tilde{z}_{i+1},$
$\ldots$, $z_{m}=\prox_{f_{m}}z_{m-1}=\prox_{f_{m}}\tilde{z}_{m-1}=\tilde{z}_{m}$,
$z_{1}=\prox_{f_{1}}z_{m}=\prox_{f_{1}}\tilde{z}_{m}=\tilde{z}_{1}$,
$z_{2}=\prox_{f_{2}}z_{1}=\prox_{f_{2}}\tilde{z}_{1}=\tilde{z}_{2}$,
$\ldots$,
$z_{i-1}=\prox_{f_{i-1}}z_{i-2}=\prox_{f_{1}}\tilde{z}_{i-2}=\tilde{z}_{i-1}.$
Hence $\bz=\tilde{\bz}$.

\ref{i:iterative}: Set $f_{m+1}=f_{1}$ and $F_{m+1}=F_{1}$. We show $F_{i+1}=\prox_{f_{i+1}}(F_{i})$ for $1\leq i\leq m$.
To this end, let $z\in F_{i}$. Then
$z=\prox_{f_{i}}\cdots\prox_{f_{1}}\prox_{f_{m}}\cdots\prox_{f_{i+1}}z,$
so
$$\prox_{f_{i+1}}z=\prox_{f_{i+1}}\prox_{f_{i}}\cdots\prox_{f_{1}}\prox_{f_{m}}\cdots\prox_{f_{i+2}}(\prox_{f_{i+1}}z)$$
which gives $\prox_{f_{i+1}}z\in F_{i+1}$. Thus, $\prox_{f_{i+1}}(F_{i})\subseteq F_{i+1}$.
To show the converse,
let $z\in F_{i+1}$. Then
\begin{equation}\label{e:cyclei}
z=\prox_{f_{i+1}}\cdots\prox_{f_{1}}\prox_{f_{m}}\cdots\prox_{f_{i+2}}z.
\end{equation}
Define $$z_{i+2}=\prox_{f_{i+2}}z,\quad z_{i+3}=\prox_{f_{i+3}}z_{i+2},\quad \ldots,
\quad z_{m}=\prox_{f_{m}}z_{m-1},$$
$$z_{1}=\prox_{f_{1}}z_{m},\quad z_{2}=\prox_{f_{2}}z_{1},\quad \ldots,
\quad z_{i}=\prox_{f_{i}}z_{i-1}.$$
Then \eqref{e:cyclei} implies $z=\prox_{f_{i+1}}z_{i}$, so that
$\bz=(z_{1},\ldots,z_{i}, z, z_{i+2},\ldots, z_{m})\in \ZZ,$
 $z_{i}\in F_{i} \text{ and  } z=\prox_{f_{i+1}}z_{i}.$
Thus, $F_{i+1}\subseteq\prox_{f_{i+1}}(F_{i})$. Hence $F_{i+1}=\prox_{f_{i+1}}(F_{i})$.
%
Altogether,
\eqref{e:image} is established.

\ref{i:twos}:  Apply \ref{i:iterative} with $m=2$.
\end{proof}

\begin{lemma}\label{l:c:clear}
 Suppose that $S=\bigcap_{i=1}^m\argmin f_{i}\neq\varnothing$. Then
\begin{enumerate}
\item \label{i:ithcycle}
$F_{i}=S$ for every $1\leq i\leq m$.
\item \label{i:truecycle}
$\ZZ=\menge{(z,\ldots,z)}{z\in S}.$
\end{enumerate}
\end{lemma}
\begin{proof} Because $\prox_{f_{i}}$ is firmly nonexpansive, $\Fix(\prox_{f_{i}})=\argmin f_{i}$, and
$S\neq\varnothing$, by \cite[Corollary 4.51]{BC2017} or \cite[Lemma 2.1]{bruck}
 for every $i$ we have
\begin{equation}\label{e:fset}
\Fix(\prox_{f_{i}}\cdots\prox_{f_{1}}\prox_{f_{m}}\cdots\prox_{f_{i+1}})=S.
\end{equation}
\ref{i:ithcycle}: This is \eqref{e:fset}.
\ref{i:truecycle}: For every $(z_{1},\ldots, z_{m})\in \ZZ$, \eqref{e:cycle} gives $z_{1}\in F_1$.
Since $F_{1}=S$ by \ref{i:ithcycle}, \eqref{e:cycle} further implies $z_{1}=z_{2}=\cdots=z_{m}$.
\end{proof}

%

Lemma~\ref{l:c:clear} shows that it is natural to assume
$S=\bigcap_{i=1}^{m}\argmin f_{i}=\varnothing$ when studying classical cycles and gap vectors.

\subsection{Classical cycles and gap vectors via the Attouch-Th\'era duality}

Using the Attouch-Th\'era duality with
$A=\partial\ff \text{ and } B=\Id-\bR,$
 and the identity $$-\Id\circ(\Id-\bR)^{-1}\circ(-\Id)=(\Id-\bR)^{-1}$$ for the
 linear relation $(\Id-\bR)^{-1}$,
we can formulate the primal-dual inclusion problem:
\begin{align}
(P) & \quad 0\in \partial\ff(\bx)+(\Id-\bR)\bx,\\
(D) & \quad 0\in (\partial\ff)^{-1}(\by)+(\Id-\bR)^{-1}\by.
\end{align}

While the primal problem ($P$) solves for the {classical cycles}, the dual ($D$) for the pair
$(\partial\ff , \Id-\bR)$ solves for the
\emph{classical gap-vectors} of $\ff$.
$(P)$ has a solution (respectively, no solution) if and only if $(D)$ has a solution (respectively, no solution).

\begin{proposition}[gap vector]\label{p:gap:v}
 The solution set of $(D)$ is at most a singleton (possibly empty).
\end{proposition}
\begin{proof}
Since
$(\Id-\bR)^{-1}=\frac{1}{2}\Id+N_{\dper}+T$ by Fact~\ref{f:shift}\ref{i:r2}
(or use Fact~\ref{l:inv:strong}),
the monotone operator
$$
(\partial\ff)^{-1}+(\Id-\bR)^{-1} =\frac{1}{2}\Id+(N_{\dper}+T+(\partial\ff)^{-1})
$$
is strongly monotone, so $[(\partial\ff)^{-1}+(\Id-\bR)^{-1}]^{-1}(0)$ is at most a singleton.
\end{proof}

\begin{proposition}[cycle and gap vectors]\label{p:classical:cg} Consider the sets of classical cycle and
 gap vectors defined respectively by
\begin{align}
\ZZ & =\menge{\bx\in\bX}{0\in \partial\ff(\bx)+(\Id-\bR)\bx},\\
\bG &=\menge{\by\in\bX}{0\in (\partial\ff)^{-1}(\by)+(\Id-\bR)^{-1}\by}.
\end{align}
We have
\begin{enumerate}
\item\label{i:cg1} $\ZZ=\bigcup_{\by\in \bG}(\Id-\bR)^{-1}(-\by)\cap (\partial\ff)^{-1}(\by)$.
\item\label{i:cg2} $\bG=\bigcup\menge{\bR \bx-\bx}{\bx\in \ZZ}.$ If $\bG\neq\varnothing$, then $\bG$ is a singleton
$\by\in\dper$ and $\by=\bR \bx-\bx$ for every $\bx\in \ZZ$.
\end{enumerate}
\end{proposition}

\begin{proof}
Apply Facts~\ref{dualityat}, \ref{f:shift}\ref{i:r2} and Proposition~\ref{p:gap:v}.
\end{proof}
%
\subsection{Relationship among fixed point sets $F_{i}$'s}
In view of Lemma~\ref{l:inter:rel}, Propositions~\ref{p:classical:cg} and \ref{p:gap:v}, we have that
$F_{i}\neq\varnothing$ if and only if $\ZZ\neq\varnothing$ if and only if $\bG\neq \varnothing$, and that
$\bG$ is either empty or a singleton. Therefore, below it is not surprising that we assume that the dual
$(D)$ has a solution.

\begin{corollary}[fixed point sets of cyclic proximal mapping] \label{c:fixset:map}
Assume that $\by=(y_{1},\ldots, y_{m})$ is the unique solution of the dual problem $(D)$.
Then the following assertions hold:
\begin{enumerate}
\item\label{i:proxmap}
For $1\leq i\leq m-1$, the mapping $\prox_{f_{i+1}}: F_{i}\rightarrow
 F_{i+1}$ is bijective and it is given by
$z\mapsto z-y_{i+1}.$
The mapping $\prox_{f_{1}}:F_{m}\rightarrow F_{1}$ is bijective and it is given by
$z\mapsto z-y_{1}.$
\item\label{i:proximg}
For $1\leq i\leq  m-1$, the fixed point sets $F_{i+1}=F_{i}-y_{i+1}$, and $F_{1}=F_{m}-y_{1}$.
Consequently, $F_{i}$ is just a translation of $F_{1}$.
\end{enumerate}
\end{corollary}
\begin{proof} By the assumption, $(D)$ has solution, hence $(P)$ has a solution by the Attouch-Th\'era duality.
In view of Lemma~\ref{l:inter:rel}, each fixed point set $F_{i}$ is nonempty.

\ref{i:proxmap}: By Lemma~\ref{l:inter:rel}\ref{i:iterative}, $\prox_{f_{i+1}}$
is onto. Suppose $z, \tilde{z}\in F_{i}$ and $\prox_{f_{i+1}}z=\prox_{f_{i+1}}\tilde{z}$.
Because
$$z=\prox_{f_{i}}\cdots\prox_{f_{1}}\prox_{f_{m}}\cdots\prox_{f_{i+1}}z,$$
$$\tilde{z}=\prox_{f_{i}}\cdots\prox_{f_{1}}\prox_{f_{m}}\cdots\prox_{f_{i+1}}\tilde{z},$$
we have $z=\tilde{z}$. Hence $\prox_{f_{i+1}}$ is injective on $F_{i}$.

Now for every $z\in F_{i}$, we have
\begin{equation}\label{e:fixedseti}
z=\prox_{f_{i}}\cdots\prox_{f_{1}}\prox_{f_{m}}\cdots\prox_{f_{i+1}}z.
\end{equation}
Set $$z_{i+1}=\prox_{f_{i+1}}z,\quad z_{i+2}=\prox_{f_{i+2}}z_{i+1}, \quad
\ldots, \quad z_{m}=\prox_{f_{m}}z_{m-1}, $$
$$z_{1}=\prox_{f_{1}}z_{m},\quad  \ldots,\quad
z_{i-1}=\prox_{f_{i-1}}z_{i-2}.$$
 Then $z=\prox_{f_{i}}z_{i-1}$ by \eqref{e:fixedseti}, and
$$\bz=(z_{1},\ldots, z_{i-1}, z,z_{i+1},\ldots, z_{m})$$
is a solution to $(P)$. By Proposition~\ref{p:classical:cg}\ref{i:cg2},
$$(z_{m}-z_{1},z_{1}-z_{2}, \ldots, z-z_{i+1}, \ldots, z_{m-1}-z_{m})=\bR \bz-\bz =\by.$$
Thus,
$z-z_{i+1}=y_{i+1}$ so that $z_{i+1}=z-y_{i+1}$, i.e.,
$\prox_{f_{i+1}}z=z-y_{i+1}$.

The proof for $\prox_{f_{1}}$ is analogous.

\ref{i:proximg}: This immediate from \ref{i:proxmap}.
\end{proof}
\begin{remark} Corollary~\ref{c:fixset:map} extends \cite[Proposition 3.2(v)]{reich} and \cite[Theorem 3.3(v)]{wangb}
 from compositions
of two resolvents
to compositions of multi-proximal mappings.
\end{remark}

\begin{corollary}\label{c:notsuzuki}
 Suppose that one of the fixed point set $F_{i}$s is nonempty and bounded, e.g., when one of
$\dom f_{i}$s is bounded. Then the following are equivalent:
\begin{enumerate}
\item\label{i:rads} $F_{1}=\cdots=F_{m}\neq\varnothing$.
\item\label{i:mini} $\bigcap_{i=1}^{m}\argmin f_{i}\neq\varnothing$.
\end{enumerate}
Under either of the two conditions, we have
\begin{equation}\label{e:suzuki:christ}
F_{1}=\cdots=F_{m}=\bigcap_{i=1}^{m}\argmin f_{i}=\bigcap_{i=1}^{m}\Fix(\prox_{f_{i}})\neq\varnothing.
\end{equation}
\end{corollary}
\begin{proof}
\ref{i:rads}$\Rightarrow$\ref{i:mini}: By Lemma~\ref{l:inter:rel}\ref{i:FZ} and Proposition~\ref{p:classical:cg},
the dual problem $(D)$ has a unique solution, say $\by$.
Because that $F_{i}=F_{i+1}=F_{i}-y_{i+1}$ by
Corollary~\ref{c:fixset:map}\ref{i:proximg},
and that $F_{i}$ is bounded, closed
and convex, Radstrom's Cancellation Theorem \cite[page 68]{BC2017} implies $y_{i+1}=0$ for $i=1,\ldots, m$,
so $\by=0$. In view of Proposition~\ref{p:classical:cg}\ref{i:cg2} with $\bx$ being the solution to \eqref{e:sub:form},
$\bR\bx-\bx=0$, thus,
$x_{1}=\cdots=x_{m}$. By \eqref{e:sub:form}, $(0,\ldots, 0)\in\partial f_{1}(x_{1})\times
\cdots\times \partial f_{m}(x_{1})$,
so $x_{1}\in\bigcap_{i=1}^{m}\argmin f_{i}$. This establishes \ref{i:mini}.

\ref{i:mini}$\Rightarrow$\ref{i:rads}: Apply Lemma~\ref{l:c:clear}\ref{i:ithcycle}.

Equation \eqref{e:suzuki:christ} follows from Lemma~\ref{l:c:clear}\ref{i:ithcycle} and
$\argmin f_{i}=\Fix\prox_{f_{i}}$.
\end{proof}

\begin{corollary}\label{c:f:notempty}
 Suppose that $\argmin f_{i}\neq\varnothing$ for $1\leq i\leq m$. Then the following are equivalent:
\begin{enumerate}
\item\label{i:f:nonempty1} $F_{1}=\cdots=F_{m}\neq\varnothing$.
\item\label{i:f:nonempty2} $\bigcap_{i=1}^{m}\argmin f_{i}\neq\varnothing$.
\end{enumerate}
Under either of the two conditions, we have
$
F_{1}=\cdots=F_{m}=\bigcap_{i=1}^{m}\argmin f_{i}=\bigcap_{i=1}^{m}\Fix(\prox_{f_{i}})\neq\varnothing.
$
\end{corollary}
\begin{proof}
\ref{i:f:nonempty1}$\Rightarrow$\ref{i:f:nonempty2}:
By Lemma~\ref{l:inter:rel}\ref{i:FZ} and Proposition~\ref{p:classical:cg},
the dual problem $(D)$ has a unique solution, say $\by$.
The mapping $\prox_{f_{i+1}}: F_{i}\rightarrow F_{i+1}$
is given by $z\rightarrow z-y_{i+1}$ by
Corollary~\ref{c:fixset:map}. Using the assumption $F_{i}=F_{i+1}$, we have
\begin{equation}\label{e:blowup}
(\forall n\in\NN)\ \prox_{f_{i+1}}^{n}(z)=z-ny_{i+1}\in F_{i}.
\end{equation}
Take $x\in \argmin f_{i+1}=\Fix\prox_{f_{i+1}}$, which is possible by the assumption.
As $\prox_{f_{i+1}}$ is nonexpansive,
$$\|\prox_{f_{i+1}}^{n}(z)-x\|=\|\prox_{f_{i+1}}^{n}(z)-\prox_{f_{i+1}}^{n}(x)\|\leq \|z-x\|,$$
so $(\prox_{f_{i+1}}^{n}(z))_{n\in\NN}$ is bounded, and this implies $y_{i+1}=0$ by \eqref{e:blowup}.
Since the same arguments apply to every
$\prox_{f_{i}}$ for $1\leq i\leq m$, we obtain $\by=0$.
In view of Proposition~\ref{p:classical:cg}\ref{i:cg2}
with $\bx$ being the solution to \eqref{e:sub:form},
$\bR\bx-\bx=0$, thus,
$x_{1}=\cdots=x_{m}$. By \eqref{e:sub:form}, $(0,\ldots, 0)\in\partial f_{1}(x_{1})\times
\cdots\times \partial f_{m}(x_{1})$,
so $x_{1}\in\bigcap_{i=1}^{m}\argmin f_{i}$. This establishes \ref{i:f:nonempty2}.

\ref{i:f:nonempty2}$\Rightarrow$\ref{i:f:nonempty1}: This and the remaining conclusion follow from Lemma~\ref{l:c:clear}\ref{i:ithcycle}.
\end{proof}

\begin{remark}\label{suzuki:s}
 In \cite{suzuki}, Suzuki called
\begin{equation}\label{e:bcondition}
\bigcap_{i=1}^{m}\Fix(\prox_{f_{i}})=F_{1}=\cdots=F_{m}
\end{equation}
as Bauschke's condition for the collection of proximal mappings $(\prox_{f_{i}})_{i=1}^m$.
\begin{enumerate}
\item When $\bigcap_{i=1}^{m}\Fix(\prox_{f_{i}})\neq\varnothing$, it is well-known that
\eqref{e:bcondition} holds ; cf. \cite[Corollary 4.51]{BC2017}.
\item When $\bigcap_{i=1}^{m}\Fix(\prox_{f_{i}})=\varnothing$, we have two cases to consider. Case 1: $\ZZ=\varnothing$.
Lemma~\ref{l:inter:rel} yields $F_{1}=\cdots=F_{m}=\varnothing$, hence \eqref{e:bcondition} holds because
of $\bigcap_{i=1}^{m}\Fix(\prox_{f_{i}})=\varnothing$. See, e.g.,
Example~\ref{e:no:cycle}.
Case 2: $\ZZ\neq\varnothing$. Lemma~\ref{l:inter:rel} yields $F_{i}\neq\varnothing$ for
$1\leq i\leq m$, hence \eqref{e:bcondition} fails because of $\bigcap_{i=1}^{m}\Fix(\prox_{f_{i}})=\varnothing$.
Requiring $F_{1}=\cdots=F_{m}$ implies
that $F_{1}=\cdots=F_{m}$ has to be unbounded by Corollary~\ref{c:notsuzuki}.
If each $\argmin f_{i}\neq\varnothing$, then $F_{i}\neq F_{j}$ for some $1\leq i, j\leq m$ by
Corollary~\ref{c:f:notempty}.
\end{enumerate}
\end{remark}


The question is what should one do if the classical cycle and gap vectors do not exist.

\section{The generalized cycle and gap vectors}\label{s:gcycle:gfixed}
In this section, using the Attouch-Th\'era duality we provide generalized cycles and gap vectors
for $\closu(\ff\Box\iota_{\Delta})$. One remarkable
feature is that the generalized cycles and gap vectors always exist, and they can recover the classical ones
for $\ff$
whenever the latter exist.

To deal with the absence of solutions
(or non-attainments of the classical cycle or gap vectors), we need to enlarge
the primal $(P)$ or dual $(D)$ so that it has a solution.

\subsection{Extending the dual approach}

Since the linear relation
$(\Id-\bR)^{-1}=\frac{1}{2}\Id+N_{\dper}+T$ by Fact~\ref{f:shift}\ref{i:r2}, and  $\partial\iota_{\Delta}^*=\partial\iota_{\dper}=N_{\dper},$
we have
\begin{align}
(\partial\ff)^{-1}+(\Id-\bR)^{-1} &=
\partial\ff^*+\frac{1}{2}\Id+T+\partial \iota_{\Delta}^*
=\partial\ff^*+\partial \iota_{\Delta}^*+\frac{1}{2}\Id+T \label{e:larger1}\\
& \subseteq \partial (\ff^*+\iota_{\Delta}^*)+\frac{1}{2}\Id+T\\
& =\frac{1}{2}\big[\Id+\big(2T+2\partial (\ff^*+\iota_{\Delta}^*)\big) \big].\label{e:larger2}
\end{align}
The enlarged dual
\begin{equation}\label{e:gdual}
(\tilde{D}) \quad 0\in \partial(\ff^*+\iota_{\Delta}^*)(\by)+\frac{1}{2}\by+T\by
\end{equation}
always has a unique solution. We call the unique $\by$ given by \eqref{e:gdual}
as the \emph{generalized gap vector}
of $\closu(\ff\Box\iota_{\Delta})$.

\begin{proposition}[existence and uniqueness of generalized gap vector]\label{p:dual:sol}
The following assertions hold:
\begin{enumerate}
\item\label{i:dtd1}The enlarged dual $(\tilde{D})$ always has a unique solution.
\item\label{i:dtd2} If $(D)$ has a solution, then it is exactly the solution of $(\tilde{D})$.
\end{enumerate}
\end{proposition}
\begin{proof}
\ref{i:dtd1}: 
Because
$\partial(\ff^*+\iota_{\Delta}^*)$ and $T$ are maximally monotone, we have
$2T+2\partial (\ff^*+\iota_{\Delta}^*)$ maximally monotone by \cite[Corollary 25.5]{BC2017}
or \cite[Corollary 32.3]{Simons2},
so Minty's theorem \cite[Theorem 21.1]{BC2017}
applies.

\ref{i:dtd2}: This follows from \eqref{e:larger1}--\eqref{e:larger2} and \ref{i:dtd1}.
\end{proof}

\subsection{Extending the primal approach}
One can also start from the primal
$(P)  \quad 0\in \partial\ff(\bu)+(\Id-\bR)\bu.$
Because $\partial \ff+(\Id-\bR)$ is already maximally monotone by
\cite{Simons2}, one cannot find enlargements of
($P$) that always have a solutions. We need to
rewrite it in an equivalent form then do enlargements.
In view of
$-(\Id-\bR)\bu\in \partial\ff(\bu)$, $-(\Id-\bR)\bu\in\dper$, Lemmas~\ref{l:exact} and \ref{l:closed},
we have
\begin{align}
-(\Id-\bR)\bu
& \in \partial\ff(\bu)\cap\dper
=\partial\ff(\bu)\cap\partial\iota_{\Delta}(\dd)\subseteq\partial(\ff\Box\iota_{\Delta})(\bu+\dd)\label{e:pr1}\\
&\subseteq\partial [\closu{(\ff\Box\iota_{\Delta})}](\bu+\dd),\label{e:pr2}
\end{align}
where $\dd\in\Delta$. As $(\Id-\bR)(\dd)=0$, we can write equations \eqref{e:pr1}
as
$$0\in \partial (\ff\Box\iota_{\Delta})(\bu+\dd)+(\Id-\bR)(\bu+\dd).$$
However, because $f\Box\iota_{\Delta}$ might not be lower semicontinuous,
$\partial (\ff\Box\iota_{\Delta})$ need not be maximally monotone. Therefore, by \eqref{e:pr2}
we
consider the enlargement
$$0\in \partial [\closu{(\ff\Box\iota_{\Delta})}](\bu+\dd)+(\Id-\bR)(\bu+\dd).$$
With $\dd=-P_{\Delta}(\bu)
\in\Delta$ and
$\bx=\bu+\dd=P_{\Delta^{\perp}}(\bu)\in \dper$, we have
\begin{equation}\label{e:pap}
0\in \partial [\closu{(\ff\Box\iota_{\Delta})}](\bx)+(\Id-\bR)(\bx), \text{ and } \bx\in \dper.
\end{equation}
The solution $\bx$ given by \eqref{e:pap} is called a \emph{generalized cycle} of $\closu(\ff\Box
\iota_{\Delta})$.
This is so-called because of \eqref{e:ab:cycle} and that \eqref{e:pap} can be reformulated as
$\bx=\prox_{\closu(\ff\Box
\iota_{\Delta})}(\bR\bx) \text{ and } \bx\in \dper.$
One amazing property of \eqref{e:pap} is that
it always has a solution, see Theorem~\ref{t:pd}\ref{i:pd1} below!

%

\subsection{The primal-dual approach}

The generalized cycles and gap vectors of $\closu(\ff\Box\iota_{\Delta})$ can be put into the framework of
the Attouch-Th\'era duality.

\begin{theorem}[generalized cycle and gap vectors via duality]\label{t:pd}
Consider the following Attouch-Th\'era primal-dual problems
\begin{align}
(\tilde{P}) & \quad 0 \in \partial[\closu(\ff\Box\iota_{\Delta})](\bx)+(\Id-\bR)\bx, \text{ and } \bx\in\dper,\label{e:p:form} \\
(\tilde{D}) &\quad 0\in \partial(\ff^*+\iota_{\Delta}^*)(\by)+\frac{1}{2}\by+T\by. \label{e:d:form}
\end{align}
Then the following assertions hold:
\begin{enumerate}
\item\label{i:pd2} ($\tilde{D}$) has a unique solution.
\item\label{i:pd1} $(\tilde{P})$ is the Attouch-Th\'era dual of $(\tilde{D})$ associated with the pair $(\partial(\ff^*+\iota_{\Delta}^*), \Id/2+T)$, and ($\tilde{P}$) has a unique solution.
\end{enumerate}
\end{theorem}
\begin{proof}
\ref{i:pd2}: ($\tilde{D}$) has a unique solution by Proposition~\ref{p:dual:sol}\ref{i:dtd1}.

\ref{i:pd1}:
Let us compute the Attouch-Th\'era dual of $(\tilde{D})$
associated with the pair $(\partial(\ff^*+\iota_{\Delta}^*), \Id/2+T)$. Because $(\ff^*+\iota_{\Delta}^*)^*=\closu(\ff\Box
\iota_{\Delta})$ and $\partial((\ff^*+\iota_{\Delta}^*)^*)=(\partial(\ff^*+\iota_{\Delta}^*))^{-1}$,
we have the Attouch-Th\'era dual
\begin{equation}\label{e:d:inverse}
\quad 0\in \partial[\closu(\ff\Box
\iota_{\Delta})](\bx)+\left(\frac{1}{2}\Id+T\right)^{-1}(\bx).
\end{equation}
Since $(\frac{1}{2}\Id+T)^{-1}=\Id-\bR+2P_{\dD}$ by Fact~\ref{f:shift}\ref{i:r3}, we obtain
$$\quad 0 \in \partial[\closu(\ff\Box
\iota_{\Delta})](\bx)+(\Id-\bR)\bx+2P_{\dD}(\bx).$$
Because
$\ran(\Id-\bR)\subseteq\dper$ and Lemma~\ref{l:range}, the above implies
$$ -2P_{\dD}(\bx) \in \partial[\closu(\ff\Box
\iota_{\Delta})](\bx)+(\Id-\bR)\bx \subseteq\dper,$$
from which $2P_{\dD}(\bx)\in \dD\cap\dper$, so $P_{\dD}(\bx)=0$, and $\bx\in\dper$.
Hence, \eqref{e:d:inverse} is equivalent to
\begin{equation}\label{e:dual}
0 \in \partial[\closu(\ff\Box
\iota_{\Delta})](\bx)+(\Id-\bR)\bx, \text{ and } \bx\in\dper,
\end{equation}
which is precisely \eqref{e:p:form}.

Now ($\tilde{P}$) has solutions because of \ref{i:pd2} and Fact~\ref{dualityat}.
To prove that the solution
of ($\tilde{P}$) is unique, let
$\bx, \tilde{\bx}$ be two solutions to ($\tilde{P}$). Since
$-(\Id-\bR)(\bx), -(\Id-\bR)(\tilde{\bx})$ have to be the unique solution to ($\tilde{D}$),
we have
$\bx-\tilde{\bx}\in\ker(\Id-\bR)=\Delta$. Because
$\bx-\tilde{\bx}\in\dper$, we conclude that $\bx-\tilde{\bx}=0$.
\end{proof}

\subsection{Relationship between the classical cycles and generalized cycles}

It is natural to ask how does the solution of ($\tilde{P}$) relate to the classical cycles of $\ff$.

\begin{theorem}\label{t:one} Let $\bx$ be the generalized cycle of $\closu(\ff\Box\iota_{\Delta})$, i.e.,
$$(\tilde{P}) \quad 0 \in \partial[\closu(\ff\Box\iota_{\Delta})](\bx)+ (\Id-\bR)\bx, \text{ and } \bx\in\dper.$$
If $\closu(\ff\Box\iota_{\Delta})(\bx)=(\ff\Box\iota_{\Delta})(\bx)$ and
$\ff\Box\iota_{\Delta}$ is exact at $\bx$,
then $\exists \bu\in\bX$ such that $\bx=\bu+\bv$, $\bv=-P_{\Delta}(\bu)\in\dD$,
and
$$0\in \partial \ff (\bu)+(\Id-\bR)\bu.$$
Consequently, $\bu$ is a classical cycle for $\ff$.
\end{theorem}
\begin{proof} As $\closu(\ff\Box\iota_{\Delta})(\bx)=(\ff\Box\iota_{\Delta})(\bx)$,
by Lemma~\ref{l:closed}
we have
\begin{equation}\label{e:remove:clos}
\partial [\closu(\ff\Box\iota_{\Delta})](\bx)=\partial (\ff\Box\iota_{\Delta})(\bx).
\end{equation}
Because $\ff\Box\iota_{\Delta}$ is exact at $\bx$,
there exist $\bu\in \bX, \bv\in\dD$ such that
$\bx=\bu+\bv$ and $(\ff\Box\iota_{\Delta})(\bx)=\ff(\bu)+\iota_{\Delta}(\bv)$.
Using \eqref{e:remove:clos} and Lemmas~\ref{l:exact}, we obtain
$$\partial[\closu(\ff\Box\iota_{\Delta})](\bx)=\partial(\ff\Box\iota_{\Delta})(\bx)
=\partial\ff(\bu)\cap N_{\dD}(\bv)=\partial\ff(\bu)\cap \dper.$$
Also
$(\Id-\bR)(\bu+\bv)=(\Id-\bR)\bu$ because $\bR\bv=\bv$ and $\bv\in\dD$. Hence, ($\tilde{P}$)
simplifies to
$$0\in \partial\ff(\bu)+(\Id-\bR)\bu.$$
Finally, $\bx=\bu+\bv$, $\bx\in\dper$ and $\bv\in\Delta$ imply
$0=P_{\Delta}(\bx)=P_{\Delta}\bu+\bv$, so $\bv=-P_{\Delta}(\bu).$
\end{proof}

\begin{theorem}\label{t:two} 
Let $\bu$ be a classical cycle for $\ff$, i.e.,
\begin{equation}\label{e:c-cycle}
0\in \partial\ff(\bu)+(\Id-\bR)\bu.
\end{equation}
Set
$\bv=-P_{\Delta}(\bu)\in\dD$
and $\bx=\bu+\bv=P_{\Delta^{\perp}}(\bu)$.
Then
\begin{enumerate}
\item\label{i:lsc:exact} $\ff\Box\iota_{\Delta}$ is lower semicontinuous and exact at $\bx$.
\item \label{i:gen:cycle}
$\bx \in \dper$ and $\bx$ solves
\begin{equation}\label{e:t:cycle}
(\tilde{P}) \quad 0 \in \partial (\ff\Box\iota_{\Delta})(\bx)+(\Id-\bR)\bx
=\partial[\closu(\ff\Box\iota_{\Delta})](\bx)+(\Id-\bR)\bx.
\end{equation}
Consequently, $\bx$ is the generalized cycle for $\closu(\ff\Box\iota_{\Delta})$.
\item\label{i:t:cycleset} $P_{\Delta^{\perp}}(\ZZ)=\{\bx\}$.
\end{enumerate}
\end{theorem}

\begin{proof}
\ref{i:lsc:exact}:
As $\bR\bu-\bu\in\dper$, \eqref{e:c-cycle} implies that $ \bR\bu-\bu\in\partial\ff(\bu)\cap\dper=
\partial\ff(\bu)\cap\partial\iota_{\Delta}(\bv)\neq\varnothing$.
By Lemma~\ref{l:exact}\ref{i:exact:sub},
$\ff\Box\iota_{\Delta}$ is exact at $\bx=\bu+\bv$, i.e.,
$(\ff\Box\iota_{\Delta})(\bx)=\ff(\bu)+\iota_{\Delta}(\bv)$. Apply Lemma~\ref{l:exact}\ref{i:exact:funct}
to
obtain $\partial(\ff\Box\iota_{\Delta})(\bx)
=\partial\ff(\bu)\cap\dper\neq\varnothing$.  Then by Lemma~\ref{l:closed} we have
$\closu(\ff\Box\iota_{\Delta})(\bx)=(\ff\Box\iota_{\Delta})(\bx)$ and
\begin{equation}\label{e:cl:sub}
\partial[\closu(\ff\Box\iota_{\Delta})](\bx)=\partial(\ff\Box\iota_{\Delta})(\bx)=\partial\ff(\bu)\cap \dper.
\end{equation}

\ref{i:gen:cycle}:
Clearly $\bx\in\dper$.
Moreover,
\begin{equation}\label{e:xu}
(\Id-\bR)\bx=(\Id-\bR)(\bu)+(\Id-\bR)(\bv)=(\Id-\bR)(\bu)\in\dper.
\end{equation}
Using \eqref{e:cl:sub} and \eqref{e:xu}, we can rewrite \eqref{e:c-cycle} as
\begin{align}
0\in \partial\ff(\bu)+(\Id-\bR)\bu &\Rightarrow 0\in \partial\ff(\bu)\cap \dper+(\Id-\bR)\bx\\
& \Rightarrow 0\in \partial(\ff\Box\iota_{\Delta})(\bx)+(\Id-\bR)\bx=
\partial[\closu(\ff\Box\iota_{\Delta})](\bx)+(\Id-\bR)\bx,
\end{align}
which is \eqref{e:t:cycle}.

\ref{i:t:cycleset}: Clear from \ref{i:gen:cycle}.
\end{proof}

The following result summarizes the relationship among the classical cycles, generalized cycles, and
generalized gap vectors.

\begin{corollary}\label{c:pd:sol} With $\bx$ and $\by$ given in Theorem~\ref{t:pd},
the following assertions hold:
\begin{enumerate}
\item\label{i:pd:sol} $\bx, \by\in\dper$.
\item\label{i:by} $\by=\bR \bx-\bx$.
\item \label{i:bx} $\bx=-\frac{\by}{2}-T\by$.
\item\label{i:z:char} $\ZZ=(\bx+\Delta)\cap (\partial\ff)^{-1}(\bR \bx-\bx)=
(\Id-\bR)^{-1}(-\by)\cap(\partial f)^{-1}(\by)$.
\item \label{i:z:sub}
$\ZZ\subseteq (F_{1}\times\cdots\times F_{m})\cap (\Id-\bR)^{-1}(-\by)$.
\end{enumerate}
\end{corollary}
\begin{proof}
\ref{i:pd:sol}: This follows from Theorem~\ref{t:pd} and $\dom (\ff^*+\iota_{\Delta}^*)\subseteq
\dom\iota_{\Delta}^*=\dom\iota_{\Delta^{\perp}}=\dper$.

\ref{i:by}: Set $v=\bR \bx-\bx$. By the assumption, $-v=\bx-\bR\bx,
v\in \partial[\closu(\ff\Box\iota_{\Delta})](\bx)$, from which
$\bx\in \partial(\ff^*+\iota_{\Delta}^*)(v)$ and $-\bx\in (\Id-\bR)^{-1}(v)$.
In view of Fact~\ref{f:shift}\ref{i:r2},
it follows that
\begin{align}
0 &\in \partial(\ff^*+\iota_{\Delta}^*)(v)+(\Id-\bR)^{-1}(v)
=\partial(\ff^*+\iota_{\Delta}^*)(v)+\frac{v}{2}+N_{\dper}(v)+ T(v)\\
&=\partial(\ff^*+\iota_{\Delta}^*)(v)+\frac{v}{2}+\Delta+ T(v)=[\partial(\ff^*+\iota_{\Delta}^*)(v)+
\Delta]+\frac{v}{2}+ T(v)\\
&=\partial(\ff^*+\iota_{\Delta}^*)(v)+\frac{v}{2}+ T(v),
\end{align}
which implies $v\in\dper$ and $v$ solves ($\tilde{D}$).
By Theorem~\ref{t:pd}\ref{i:pd2}, we conclude $v=\by$.

\ref{i:bx}: Set $v=-\frac{\by}{2}-T\by$. Then
$-\by\in (\frac{\Id}{2}+T)^{-1}(v)$ and, by the assumption,
$$\by\in [\partial(\ff^*+\iota_{\Delta}^*)]^{-1}(v)=\partial[\closu(\ff\Box\iota_{\Delta})](v).$$
In view of Fact~\ref{f:shift}\ref{i:r3}, it follows that
\begin{align}
0& \in \partial[\closu(\ff\Box\iota_{\Delta})](v)+\left(\frac{\Id}{2}+T\right)^{-1}(v)\\
&=\partial[\closu(\ff\Box\iota_{\Delta})](v)+(\Id-\bR+2P_{\Delta})(v).
\end{align}
Then
\begin{equation}\label{e:bx:form}
-2P_{\Delta}(v)\in \partial[\closu(\ff\Box\iota_{\Delta})](v)+(\Id-\bR)v\subseteq \dper
\end{equation}
gives $2P_{\Delta}(v)\in\Delta\cap\dper$, so $P_{\Delta}(v)=0$ and $v\in\dper$. Therefore,
\eqref{e:bx:form} reduces to
$$0\in \partial[\closu(\ff\Box\iota_{\Delta})](v)+(\Id-\bR)v, \text{ and } v\in\dper.$$
By Theorem~\ref{t:pd}\ref{i:pd1}, $v=\bx$.


\ref{i:z:char}:
Suppose $u\in \ZZ$. By Theorem~\ref{t:two}, $\bx=\bu+\bv$ solves \eqref{e:t:cycle}, $\bv\in\Delta$,
and $\partial (\ff\Box\iota_{\Delta})(\bx)=\partial\ff(\bu)\cap\dper$.
We have
$\bu=\bx-\bv\in \bx+\Delta$ and
$$0\in \partial(\ff\Box\iota_{\Delta})(\bx)+\bx-\bR\bx\quad \Rightarrow\quad 0\in \partial\ff(\bu)+\bx-\bR\bx,$$
which implies $\bu\in (\partial\ff)^{-1}(\bR\bx-\bx)$. Then $\bu\in (\bx+\Delta)\cap (\partial\ff)^{-1}(\bR\bx-\bx).$

Conversely, let $\bu\in (\bx+\Delta)\cap (\partial\ff)^{-1}(\bR \bx-\bx)$. We have
$\bu\in \bx-\Delta$ and $\bu\in (\partial\ff)^{-1}(\bR\bx-\bx)$. Then
$\bR\bu-\bu=\bR\bx-\bx\in\partial \ff(\bu)$ so that
$0\in \partial \ff(\bu)+ \bu-\bR\bu.$
Hence $\bu\in \ZZ$.

The second equality follows from \ref{i:by} and
$(\Id-\bR)^{-1}(-\by)=\bx+\Delta.$

\ref{i:z:sub}: Apply \ref{i:z:char} and Lemma~\ref{l:inter:rel}\ref{i:proj}.
\end{proof}
\begin{remark} Corollary~\ref{c:pd:sol}\ref{i:by}\&\ref{i:bx} can also be proved by
using \cite[Remark 5.4]{bot} via paramonotonicity.
\end{remark}

Recall that the parallel sum of $\partial f$ and $\partial g$ is defined by  $(\forall x\in X)\
(\partial f\Box\partial g)(x)=\bigcup_{x=u+v}\partial f(u)\cap\partial g(v)$.
\begin{corollary}[existence of classical cycles]\label{c:znonempty}
 Let $\bx$ be given in Theorem~\ref{t:pd}. Then the following are equivalent:
\begin{enumerate}
\item\label{i:zz}
 $\ZZ\neq\varnothing$.
\item\label{i:fboxi} $[\closu(\ff\Box\iota_{\Delta})](\bx)=(\ff\Box\iota_{\Delta})(\bx)$
and $\ff\Box\iota_{\Delta}$ is exact at $\bx$.
\item\label{i:parallels}
 $(\partial \ff\Box\partial\iota_{\Delta})(\bx)\neq \varnothing$.
\end{enumerate}
\end{corollary}
\begin{proof} Observe that
\begin{equation}\label{e:dualsol}
\bx \text{ solving ($\tilde{P}$) implies }
\partial [\closu(\ff\Box\iota_{\Delta})](\bx)\neq \varnothing.
\end{equation}

\ref{i:zz}$\Leftrightarrow$\ref{i:fboxi}: Combine Theorems~\ref{t:one} and \ref{t:two}.

\ref{i:fboxi}$\Rightarrow$\ref{i:parallels}:  Lemma~\ref{l:closed} gives
\begin{equation}\label{e:noclosure}
\partial (\ff\Box\iota_{\Delta})(\bx)=\partial [\closu(\ff\Box\iota_{\Delta})](\bx).
\end{equation}
Since $\ff\Box\iota_{\Delta}$ is exact at $\bx$, by Lemma~\ref{l:exact}\ref{i:exact:funct}
there exist $\bu, \bv\in\bX$ such that
$\bx=\bu+\bv$ and
\begin{equation}\label{e:psum}
\partial (\ff\Box\iota_{\Delta})(\bx)=\partial\ff(\bu)\cap\partial\iota_{\Delta}(\bv).
\end{equation}
Combining \eqref{e:dualsol}, \eqref{e:noclosure} and \eqref{e:psum},  we get
$(\partial\ff\Box\partial\iota_{\Delta})(\bx)\neq \varnothing.$

\ref{i:parallels}$\Rightarrow$\ref{i:fboxi}: $(\partial \ff\Box\partial\iota_{\Delta})(\bx)\neq \varnothing$ implies
there exist $\bu, \bv\in\bX$ such that $\bx=\bu+\bv$ and
$\partial\ff(\bu)\cap\partial\iota_{\Delta}(\bv)\neq\varnothing$. By Lemma~\ref{l:exact}\ref{i:exact:sub},
$\ff\Box\iota_{\Delta}$ is exact at $\bx$, and $\partial (\ff\Box\iota_{\Delta})(\bx)=\partial\ff(\bu)\cap\partial\iota_{\Delta}(\bv)
\neq\varnothing$, so $\ff\Box\iota_{\Delta}$ is lower semicontinuous at $\bx$
by Lemma~\ref{l:closed}\ref{i:nonemptys}.
\end{proof}


\subsection{More on generalized cycles}
Let us reconsider
($\tilde{P}$) without the restriction $\bx\in\dper$, namely,
\begin{equation}\label{e:no:restrict}
 0 \in \partial[\closu(\ff\Box\iota_{\Delta})](\bx)+(\Id-\bR)\bx.
 \end{equation}
 Denote the set of solutions to \eqref{e:no:restrict} by $E$.
 In order to study the set $E$, we shall
 need one lemma.
\begin{lemma}\label{l:new:dual}
 The Attouch-Th\'era dual of \eqref{e:no:restrict} for the pair
 $(\partial[\closu(\ff\Box\iota_{\Delta})], \Id-\bR)$ is
\begin{equation}\label{e:new:dual}
(\tilde{D}) \quad 0\in \partial(\ff^*+\iota_{\Delta}^*)(\by)+\frac{1}{2}\by+T\by.
\end{equation}
In particular, ($\tilde{P}$) and \eqref{e:no:restrict} have exactly the same dual problem.
\end{lemma}
\begin{proof} The Attouch-Th\'era dual of \eqref{e:no:restrict} for the pair
 $(\partial[\closu(\ff\Box\iota_{\Delta})], \Id-\bR)$  is
\begin{equation}\label{e:dsolve}
0\in \partial(\ff^*+\iota_{\Delta}^*)(\by)+(\Id-\bR)^{-1}(\by).
\end{equation}
Now
$(\Id-\bR)^{-1}=\frac{\Id}{2}+T+N_{\dper}$ by Fact~\ref{f:shift}\ref{i:r2} and
$\iota_{\Delta}^*=\iota_{\dper}$,
\eqref{e:dsolve} becomes
\begin{equation}\label{e:d:exp}
0\in \partial(\ff^*+\iota_{\dper})(\by)+\frac{\by}{2}+T\by+N_{\dper}(\by)=\partial(\ff^*+\iota_{\dper})(\by)+N_{\dper}(\by)+
\frac{\by}{2}+T\by.
\end{equation}
Because $\dom (\ff^*+\iota_{\dper})=\dom\ff^*\cap\dper \subseteq\dper$, $N_{\dper}(\by)=\Delta$ if $\by\in\dper$ and
$\varnothing$ if $\by\not\in\dper$, and $0\in\Delta$, we have
\begin{equation}\label{e:s:sub}
\partial(\ff^*+\iota_{\dper})(\by)+N_{\dper}(\by)=\partial(\ff^*+\iota_{\dper})(\by).
\end{equation}
Hence \eqref{e:new:dual} holds by combining \eqref{e:d:exp} and \eqref{e:s:sub}.
\end{proof}

\begin{theorem}\label{t:many} Let $\bx$ be given in Theorem~\ref{t:pd}. Then
\begin{enumerate}
\item \label{i:E1}
$E=\bx +\Delta$. Consequently, the problem \eqref{e:no:restrict} always has infinitely many solutions.
\item\label{i:E2} $P_{\Delta^{\perp}}(E)=\{\bx\}$.
\end{enumerate}
\end{theorem}
\begin{proof}
\ref{i:E1}: Let $\dd\in\Delta$. Since
$(\Id-\bR)(\bx+\dd)=(\Id-\bR)(\bx),$
Lemma~\ref{c:subdiff} gives that $\bx+\dd\in E$. Then $\bx+\Delta\subseteq E$.
To show that $\bx+\Delta\supseteq E$, let $\tilde{\bx}\in E$. By Theorem~\ref{t:pd} and Lemma~\ref{l:new:dual},
the uniqueness of solution to ($\tilde{D}$) implies $\by=\bR\bx-\bx=\bR\tilde{\bx}-\tilde{\bx}$, from
which $\tilde{\bx}-\bx\in \ker(\Id-\bR)=\Delta$, so
$\tilde{\bx}\in\bx+\Delta$. Then $E\subseteq\bx+\Delta$. Altogether $E=\bx+\Delta.$

\ref{i:E2}: Use \ref{i:E1} and $\bx\in\dper$.
\end{proof}


\section{When is $\ff\Box\iota_{\Delta}$ lower semicontinuous and exact?}\label{s:lowerh:exact}
When $\ff\Box\iota_{\Delta}$ is closed and exact, we have that
the set of classical cycles of $\ff$ is nonempty by Corollary~\ref{c:znonempty}, and  that
Theorem~\ref{t:one} allows us to find the classical
cycle of $\ff$. Thus, it is important to give some conditions under which
$\ff\Box\iota_{\Delta}$ is closed and exact.

\begin{proposition}\label{p:infexact0}
Suppose that one of the following holds:
\begin{enumerate}
\item\label{i:1} The conical hull of $(\dom f_{1}^*\times\cdots\times
\dom f_{m}^*)-\dper  \text{ is a closed subspace.}$
\item\label{i:2} The function $d\mapsto
(f_{1}(d)+\cdots+f_{m}(d))$ is coercive and the canonical hull of $$(\dom f_{1}\times\cdots\times\dom f_{m})-\Delta
\text{ is a closed subspace.} $$
\item\label{i:4} $\dom f_{1}^*\times\cdots\times\dom f_{m}^*=\bX$.
\item\label{i:6} $f_{1}^*\oplus\cdots\oplus f_{m}^*$ is continuous at some point in $\dper.$
\item\label{i:7} $X$ is finite-dimensional and $(\reli \dom f_{1}^*\times\cdots\times\reli\dom f_{m}^*)\cap \dper\neq\varnothing$.
\end{enumerate}
Then
$\ff\Box\iota_{\Delta}$ is proper, lower semicontinuous and convex, and exact on $\bX$.
\end{proposition}
\begin{proof}
Apply \cite[Proposition 15.7]{BC2017} and \cite[Proposition 15.5]{BC2017} to $\ff$ and $\iota_{\Delta}$,
and use $\ff^*=f_{1}^*\oplus\cdots\oplus f_{m}^*$ and $\iota_{\Delta}^*=\iota_{\dper}.$
\end{proof}

\begin{proposition}\label{p:infexact}
 Let $\inf f_{i}>-\infty$ for $i=1,\ldots, m$, and suppose that
at least one of $f_{i}$ is coercive. Then
$\ff\Box\iota_{\Delta}$ is proper, lower semicontinuous and convex, and exact on $\bX$.
\end{proposition}
\begin{proof} Without loss of any generality, we can assume that $f_{1}$ is coercive. Because
$\dom (\ff\Box\iota_{\Delta})=\dom\ff+\Delta$ and $\inf f_{i}>-\infty$
for all $i$, the function $\ff\Box\iota_{\Delta}$ is proper, and convex by
\cite[Proposition 12.11]{BC2017}.

First, we show that $\ff\Box\iota_{\Delta}$ is exact. By the definition, for $\bx=(x_{1},\ldots, x_{m})\in \bX$,
$$(\ff\Box\iota_{\Delta})(\bx)=\inf_{d\in X}\big(f_{1}(x_{1}-d)+\cdots+f_{m}(x_{m}-d)\big)=\inf_{d\in X}
g(d), $$
where $g(d)=f_{1}(x_{1}-d)+\cdots+f_{m}(x_{m}-d)$. Since $f_{1}$ is coercive, and
$$g(d) \geq f_{1}(x_{1}-d)+\sum_{i=2}^{m}\inf f_{i}$$
we have that $g$ is coercive, so $g$ has a minimizer over $X$. Hence $\ff\Box\iota_{\Delta}$ is exact.

Next, we show that $\ff\Box\iota_{\Delta}$ is lower semicontinuous. Let
$\bx=(x_{1},\ldots, x_{m})\in \bX$ and
$\bx_{n}\rightarrow\bx$ with $\bx_{n}=(x^{(n)}_{1},\ldots, x^{(n)}_{m})\in\bX$. We show that
\begin{equation}\label{e:lsc}
(\ff\Box\iota_{\Delta})(\bx)\leq\liminf_{n\rightarrow\infty} (\ff\Box\iota_{\Delta})(\bx_{n}).
\end{equation}
Set $\mu=\liminf_{n\rightarrow\infty} (\ff\Box\iota_{\Delta})(\bx_{n})$. When $\mu=+\infty$, \eqref{e:lsc}
clearly holds, so we assume that $\mu<+\infty$.
After passing to a subsequence and relabelling, we can assume that
$(\ff\Box\iota_{\Delta})(\bx_{n})\rightarrow\mu\in [-\infty, +\infty[.$
Let $(d_{n})_{n\in\NN}$ be a sequence in $X$ such that
$$(\ff\Box\iota_{\Delta})(\bx_{n})=f_{1}(x^{(n)}_{1}-d_{n})+\cdots+f_{m}(x^{(n)}_{m}-d_{n}).$$
We show that $(d_{n})_{n\in\NN}$ is bounded. Suppose this is not the case. Taking a subsequence if necessary,
we can assume that
$\|d_{n}\|\rightarrow\infty$. Then
\begin{align}
+\infty>\mu &\leftarrow (\ff\Box\iota_{\Delta})(\bx_{n})
= f_{1}(x^{(n)}_{1}-d_{n})+\cdots+f_{m}(x^{(n)}_{m}-d_{n})\\
&\geq f_{1}(x^{(n)}_{1}-d_{n})+\sum_{i=2}^{m}\inf f_{i}\rightarrow +\infty,
\end{align}
because $(\|x^{n}_{1}\|)_{n\in\NN}$ is bounded, $\|d_{n}\|\rightarrow\infty$, and $f_{1}$ is coercive.
This contradiction shows that $(d_n)_{n\in\NN}$ has to be bounded.
After passing to a subsequence and relabelling, we can
assume that $d_{n}\rightharpoonup d\in X$, i.e., $d_{n}$ converges to $d$ weakly. Then
$x^{(n)}_{i}-d_{n}\rightharpoonup x_{i}-d$,  and using weakly lower semicontinuity of $f_{i}$ we have
\begin{align}
\mu &=\lim_{n\rightarrow\infty} (\ff\Box\iota_{\Delta})(\bx_{n})
=\lim_{n\rightarrow\infty}(f_{1}(x_{1}^{(n)}-d_{n})+\cdots+f_{m}(x^{(n)}_{m}-d_{n}))\\
&\geq \sum_{i=1}^{m}\liminf_{n\rightarrow\infty}f_{i}(x_{i}^{(n)}-d_{n})
\geq \sum_{i=1}^{m}f_{i}(x_{i}-d)\geq (\ff\Box\iota_{\Delta})(\bx),
\end{align}
which gives \eqref{e:lsc}.
\end{proof}
\begin{corollary}
Suppose that $\ff$ is coercive. Then
$\ff\Box\iota_{\Delta}$ is proper, lower semicontinuous and convex, and exact.
\end{corollary}
\begin{proof} Because
$\ff=f_{1}\oplus\cdots\oplus f_{m}$ on $\bX$, we have that
$\ff$ is coercive if and only if each $f_{i}$ is coercive on $X$; in particular,
$0\in \inte\dom f_{1}^*\times\cdots\times \inte\dom f_{m}^*$, see, e.g.,
\cite[Proposition 14.16]{BC2017}.
Apply Proposition~\ref{p:infexact0}\ref{i:1} or Proposition~\ref{p:infexact}.
\end{proof}


Without the coercivity of $f_{i}$, one can exploit the structure of $f_{i}$.
Recall that when $X$ is finite-dimensional, we call $f:X\rightarrow\RX$
a \emph{polyhedral convex function} if
$f(x)=h(x)+\iota_{C}(x)$
where for every $x\in X$, $$h(x)=\max\menge{\scal{x}{b_{i}}-\beta_{i}}{i=1,\ldots, k},$$
$$C=\menge{x\in X}{\scal{x}{b_{i}}\leq\beta_{i}, i=k+1,\ldots, m}$$
with fixed $b_{i}\in X, \beta_{i}\in\RR$; see, e.g., \cite[Section 19, page 172]{Rock70}.

\begin{proposition}\label{p:polyh}
 Suppose that $X$ is finite-dimensional, and each $f_{i}:X\rightarrow\RX$ is a proper polyhedral convex
function for $i=1,\ldots, m$. If
$\ff\Box\iota_{\Delta}$ is proper, then
$\ff\Box\iota_{\Delta}$ is a polyhedral convex function (so lower semicontinuous) and
$\ff\Box\iota_{\Delta}$ is exact.
\end{proposition}
\begin{proof} By \cite[Theorem 19.4]{Rock70}, $\ff=f_{1}\oplus\cdots\oplus f_{m}$
is a proper polyhedral convex function on $\bX$. Also
$\iota_{\Delta}$ is a proper polyhedral convex function.
Because $\ff\Box\iota_{\Delta}$ is proper, \cite[Corollary 19.3.4]{Rock70}
gives the result.
\end{proof}
\begin{corollary} Suppose that $X$ is finite-dimensional, and $f_{i}=\iota_{C_{i}}$ with each
$C_{i}\subseteq X$ being a polyhedral convex set for $i=1,\ldots, m$.
Then $\ff\Box\iota_{\Delta}=\iota_{\bC+\Delta}$ is lower semicontinuous and exact; equivalently,
$\bC+\Delta$ is a closed set.
\end{corollary}
\begin{proof} Because $\iota_{\bC}\Box\iota_{\Delta}=\iota_{\bC+\Delta}$ is proper,
Proposition~\ref{p:polyh} applies.
\end{proof}

\section{Examples}\label{s:example}
In this section, we apply the results of Section~\ref{s:gcycle:gfixed} to concrete examples.
The first example illustrates the concepts of generalized cycles and gap vectors
when the classical ones do not exist.
The
second example characterizes when the set of cycles
is a singleton or infinite for a finite number of lines in a Hilbert space.
For $C_{i}\subseteq X$, $i=1, \ldots, m$, we let
$\bC=C_{1}\times \cdots\times C_{m}\subseteq \bX.$
\begin{example}\label{e:no:cycle} Let $\alpha\geq 0$. Consider
$$C_{1}=\epi\exp=\menge{(x,r)}{r\geq \exp (-x)+\alpha \text{ and } x\in \RR}, \text{ and }
C_{2}=\RR\times\{0\}.$$
Then the following assertions hold:
\begin{enumerate}
\item\label{i:no:gap} $\iota_{\bC}$  has neither a cycle nor a gap vector.
\item\label{i:yes:gap}
 $\closu (\iota_{\bC}\Box \iota_{\Delta})
=\iota_{\overline{\bC+\Delta}}$ has both a generalized cycle and a generalized
gap vector, namely $\bx=(0,\alpha/2,0,-\alpha/2),
\by=(0,-\alpha,0,\alpha)\in\RR^4$.
\end{enumerate}
\end{example}

\begin{proof} Note that $\bC=C_{1}\times C_{2}, \Delta=\menge{(u,v,u,v)}{u, v\in\RR}$
are in $\RR^2\times\RR^2=\RR^4$.

\ref{i:no:gap}: With $z_{1},z_{2}\in\RR^2$, we need to show that
\begin{equation}\label{e:noexit}
(0,0)\in \partial \iota_{C_{1}}(z_{1})\times\partial_{C_{2}}(z_{2})+(\Id-\bR)(z_{1},z_{2})
\end{equation}
has no solution. Suppose \eqref{e:noexit} has a solution. Then there exist $z_{i}\in C_{i}$ and
$z_{2}-z_{1}\in N_{C_{1}}(z_{1}), z_{1}-z_{2}\in N_{C_{2}}(z_{2})=\{0\}\times\RR$,
so $z_{2}-z_{1}=(0,y)\in N_{C_{1}}(z_{1})$ for some $y\in \RR$. Because
$$N_{C_{1}}(x,\exp(-x))=\bigcup_{\lambda\geq 0} \lambda (-\exp(-x),-1),$$
and $N_{C_{1}}(z_{1})=\{0\}$ if $z_{1}\in\inte C_{1}$,
this implies $y=0$ so that
$z_{1}=z_{2}$. This contradicts $C_{1}\cap C_{2}=\varnothing$.

Because \eqref{e:noexit} has no solution, so $\iota_{\bC}$ has no gap vector.

\ref{i:yes:gap}: $\closu (\iota_{\bC}\Box \iota_{\Delta})$ has a generalized
cycle. We show that
\begin{equation}\label{e:yes:exit}
(0,0)\in N_{\overline{\bC+\Delta}}(\bx) + (\Id-\bR)(\bx)\text{ and } \bx \in\dper
\end{equation}
has a solution $(0,\alpha/2,0,-\alpha/2)$. Indeed,
since
$$(0, \exp (-x)+\alpha/2, 0, -\alpha/2)=(x, \exp (-x)+\alpha, x, 0)+(-x, -\alpha/2, -x, -\alpha/2)\in\bC+\Delta,$$
by sending $x\rightarrow\infty$ we have $(0,\alpha/2,0,-\alpha/2)\in \overline{\bC+\Delta}\cap\dper$.
With $z_1=(0,\alpha/2)=-z_{2}$, we show that $\bx=(z_{1},z_{2})$ solves \eqref{e:yes:exit}, i.e.,
$$(0,0)\in N_{\overline{\bC+\Delta}}(0,\alpha/2, 0, -\alpha/2)+(0,\alpha, 0, -\alpha).$$
This holds because $(0,\alpha, 0,-\alpha)\in\dper$, and for $r\geq \exp(-x)+\alpha$, $x, \tilde{x}\in\RR$,
$(u,v,u,v)\in\Delta$,
\begin{align*}
& \scal{(0,\alpha, 0, -\alpha)}{\big((x, r,\tilde{x},0)+(u,v,u,v)\big)-(0,\alpha/2,0,-\alpha/2)}\\
&=\scal{(0,\alpha, 0, -\alpha)}{(x, r,\tilde{x},0)-(0,\alpha/2,0,-\alpha/2)}\\
& =\alpha(r-\alpha/2)-\alpha(\alpha/2)=\alpha(r-\alpha)\geq 0.
\end{align*}
The gap vector is given by
$\by=\bR\bx-\bx=(0,-\alpha,0,\alpha).$
\end{proof}

\begin{example}\label{e:many:cycle} Given $m$ sets in $X$:
$C_{i}=\menge{a_{i}+t_{i} b_{i}}{t_{i}\in\RR}$
where $a_{i}\in X$ and $b_{i}\in X\setminus\{0\}$ for $i=1,\ldots, m$. Then the following assertions hold:
\begin{enumerate}
\item\label{i:exist:c} $\iota_{\bC}$ always has a classical cycle, i.e., $\ZZ\neq\varnothing$.
\item \label{i:notp}
$\iota_{\bC}$ has a unique classical cycle if and only if the set of vectors $\menge{b_{i}}{i=1,\ldots, m}$
is not parallel.
\item \label{i:yesp} $\iota_{\bC}$ has infinitely many classical cycles if and only if
the set of vectors $\menge{b_{i}}{i=1,\ldots, m}$ is parallel.
\end{enumerate}
\end{example}
\begin{proof}
\ref{i:exist:c}: Using the special structures of $\bC$ and $\Delta$, one can show that $\bC+\Delta$ is closed in $\bX$.
Then $\iota_{\bC}\Box\iota_{\Delta}=\iota_{\bC+\Delta}$ is lower semicontinuous and exact, so
$\iota_{\bC}$ has a classical cycle by Corollary~\ref{c:znonempty}.

\ref{i:notp}: Write $\bu=(u_{1},\ldots, u_{m})$. $\bu$ is a cycle of $\iota_{\bC}$  means
$0\in \bu-\bR\bu+\partial\iota_{\bC}(\bu)$, i.e.,
$$-(u_{1}-u_{m})\in\partial_{C_{1}}(u_{1}), \text{ and } u_{i}-u_{i-1}\in\partial\iota_{C_{i}}(u_{i}) \text{ for }
i=2,\ldots, m.$$
As $\partial\iota_{C_{i}}(u_{i})=\menge{x\in X}{\scal{b_{i}}{x}=0}$ for $u_{i}\in C_{i}$, we have
$$\scal{b_{1}}{u_{1}-u_{m}}=0, \scal{b_{i}}{u_{i}-u_{i-1}}=0 \text{ for } i=2,\ldots, m.$$
As $u_{i}=a_{i}+t_{i}b_{i}$, the above become
\begin{align}
\scal{b_{1}}{a_{1}-a_{m}}& +\scal{b_{1}}{b_{1}}t_{1}-\scal{b_{1}}{b_{m}}t_{m} =0\label{e:s1}\\
& \quad \vdots\nonumber\\
\scal{b_{i}}{a_{i}-a_{i-1}}& -\scal{b_i}{b_{i-1}}t_{i-1}+\scal{b_{i}}{b_{i}}t_{i}=0\\
 &\quad \vdots\nonumber\\
\scal{b_{m}}{a_{m}-a_{m-1}}& -\scal{b_m}{b_{m-1}}t_{m-1}+\scal{b_{m}}{b_{m}}t_{m}=0.\label{e:s2}
\end{align}
In the form of $\bA\bt=\bb$ with
$\bt=(t_{1},\ldots, t_{m})^{\intercal}$ and
$\bb=(\scal{b_{1}}{a_{1}-a_{m}}, \scal{b_{2}}{a_{2}-a_{1}},\ldots, \scal{b_{m}}{a_{m}-a_{m-1}})^{\intercal}$,
the linear system of equations \eqref{e:s1}--\eqref{e:s2}
has
$\bA=$
$$\begin{pmatrix}
\scal{b_{1}}{b_{1}} & 0 & 0 &0 &  \cdots & 0 & -\scal{b_{1}}{b_{m}}\\
-\scal{b_{2}}{b_{1}} &\scal{b_{2}}{b_{2}} & 0 & 0 & \cdots & 0& 0\\
0 & -\scal{b_{3}}{b_{2}} &\scal{b_{3}}{b_{3}} & 0 & \cdots & 0 & 0\\
0 & 0 &\ddots &\ddots & \cdots &0 & 0 \\
0 & 0 &\cdots &-\scal{b_{i}}{b_{i-1}} & \scal{b_{i}}{b_{i}}  & 0 & 0 \\
\vdots & \vdots & \vdots & \ddots &\ddots & \ddots & 0 \\
0 & 0 & 0 & \cdots & 0 & -\scal{b_{m}}{b_{m-1}} &\scal{b_{m}}{b_{m}}
\end{pmatrix},
$$
and is consistent by \ref{i:exist:c}.
Using co-factor expansions, e.g. \cite[p. 478]{meyer2000}, we obtain the determinant
$$\mbox{det}\bA=\prod_{i=1}^{m}\scal{b_{i}}{b_{i}}-\scal{b_{1}}{b_{m}}\prod_{i=2}^{m}\scal{b_{i}}{b_{i-1}}=
\prod_{i=1}^{m}\scal{b_{i}}{b_{i}}\bigg(1-\frac{\scal{b_{1}}{b_{m}}}{\|b_{1}\|\|b_{m}\|}
\prod_{i=2}^{m}\frac{\scal{b_{i}}{b_{i-1}}}{\|b_{i}\|\|b_{i-1}\|}\bigg).$$
Then $\mbox{det}\bA=0$ if and only if
$\{b_{i}/\|b_{i}\|| i=1, \ldots, m\}$ is a collection of parallel vectors. Because
$\bA\bt=\bb$ has a unique solution if and only if $\mbox{det}\bA\neq 0$, the result
follows.

\ref{i:yesp}: This follows from \ref{i:notp}.
\end{proof}
%


\section{Finding the generalized cycle and gap vectors}\label{s:compute}
The objective of this final section is to show that one can compute the generalized cycle and gap vectors
by the forward-backward algorithms.

To solve
$0\in \partial\ff(\bx)+(\Id-\bR)(\bx),$
we start with a fixed point reformulation.
\begin{lemma}\label{l:fix:char} Let $\gamma\in\, ]0,1[$. Then
\begin{equation}\label{e:resolv}
\ZZ=\Fix\prox_{\ff}\bR=\Fix\prox_{\gamma\ff}\big((1-\gamma)\Id+\gamma\bR\big).
\end{equation}
In particular, for $\gamma =1/2$, we have
$\ZZ=\Fix\prox_{\ff/2}\big((\Id+\bR)/2\big)$;
consequently, $\ZZ$ is the fixed point set of a product of firmly nonexpansive mappings or resolvents.
\end{lemma}
\begin{proof}
We have
\begin{align}
\bx=\prox_{\ff}\bR\bx & \Leftrightarrow 0\in \partial\ff(\bx)+\bx-\bR\bx
 \Leftrightarrow 0\in \gamma\partial\ff(\bx)+\gamma(\bx-\bR\bx)\\
& \Leftrightarrow -\gamma\bx+\gamma\bR\bx\in \partial(\gamma \ff)(\bx)
\Leftrightarrow
(1-\gamma)\bx+\gamma\bR\bx\in \bx+\partial(\gamma \ff)(\bx)\\
& \Leftrightarrow \bx=\prox_{\gamma\ff}\bigg((1-\gamma)\bx+\gamma\bR\bx\bigg).
\end{align}
\end{proof}

The fixed point characterization of $\ZZ$, i.e., \eqref{e:resolv}, suggests finding $\bx\in\ZZ$ by the forward-backward
iteration
$$\bx_{n+1}=\prox_{\gamma\ff}\bigg((1-\gamma)\bx_{n}+\gamma\bR\bx_{n}\bigg)$$
where $\gamma\in ]0,1[$.
When $\ZZ\neq\varnothing$, $\bx_{n}\rightharpoonup\bx$ for some $\bx\in \ZZ$ and $\bR \bx_{n}-\bx_{n}\rightarrow\by$
the gap vector, see, e.g., \cite[Theorem 26.14]{BC2017}. When $\ZZ=\varnothing$,
$\|\bx_{n}\|\rightarrow +\infty$, see,
e.g., \cite[Fact 2.2]{BCR}, \cite[Corollary 1.4]{bruck}, \cite[page 121]{reich3},
\cite[Corollary 2.2]{bbr78},
but the asymptotic behavior of
$(\bR\bx_{n}-\bx_{n})_{n\in\NN}$ is not clear. See also
\cite[Theorem 5.23]{BC2017} for the convergence behavior of iterates of compositions of averaged operators.

In view of the possibility of $\ZZ=\varnothing$, one can consider
 the extended Attouch-Th\'era primal-dual:
\begin{align}
(\text{EP})\quad 0 & \in \partial[\closu(\ff\Box\iota_{\Delta})](\bx)+(\Id-\bR)\bx,\\
(\text{ED})\quad  0 &\in \partial(\ff^*+\iota_{\Delta}^*)(\by)+\frac{1}{2}\by+T\by.
\end{align}
Both (EP) and (ED) always have solutions. While (EP) gives all generalized cycles,
(ED) gives the unique generalized gap vector for $\closu(\ff\Box\iota_{\Delta})$.
To make the notation simple in the following proof, let us
write $$\bg=\closu(\ff\Box\iota_{\Delta}).$$
The following forward-backward iteration scheme allows us to find the extended cycles and gap
vectors. The elegance of our approach is that $(\bR\bx_{n}-\bx_{n})_{n\in\NN}$ always
converge strongly to the unique extended gap vector, which is precisely the classical one if the latter
exists; see Proposition~\ref{p:dual:sol}\ref{i:dtd2}.
\begin{theorem}\label{t:solve} Let $\gamma\in\, ]0, 1[$, $\delta=2-\gamma$, and
let $(\lambda_{n})_{n\in\NN}$
be a sequence in $]0, \delta[$ such that
$\sum_{n\in\NN}\lambda_{n}(\delta-\lambda_{n})=+\infty$. Let
$\bx_{0}\in \bX$ and set
\begin{equation}\label{e:fwbw}
\begin{array}{l}
\text{for}\;n=0,1,\ldots\\
\left\lfloor
\begin{array}{l}
\by_{n}=(1-\gamma)\bx_{n}+\gamma\bR \bx_{n},\\
\bx_{n+1}=\bx_{n}+\lambda_{n}(\prox_{\gamma \bg}\by_{n}-\bx_{n}).
\end{array}
\right.\\[2mm]
\end{array}
\end{equation}
Then the following assertions hold:
\begin{enumerate}
\item\label{i:findg1}
 $(\bx_{n})_{n\in\NN}$ converges weakly to $\tilde{\bx}$, a generalized cycle of $\bg$, i.e., a solution of (EP).
\item\label{i:findg2}
$(\bR x_{n}-\bx_{n})_{n\in\NN}$ converges strongly to $\by=\bR\tilde{\bx}-\tilde{\bx}$, the unique generalized
gap vector of $\bg$,
i.e., the solution of (ED).
\item\label{i:findg3}
$(P_{\Delta^{\perp}}\bx_{n})_{n\in\NN}$ converges strongly to
$P_{\Delta^{\perp}}\tilde{\bx}=\bx$, the unique solution to $(\tilde{P})$.
\end{enumerate}
\end{theorem}
\begin{proof} \ref{i:findg1}\&\ref{i:findg2}:
Let $A=\partial \bg$ and $B=\Id-\bR$. Then
$A:\bX\To\bX$ is maximally monotone, and
$B:\bX\rightarrow\bX$ is $1/2$-cocercive by Fact~\ref{l:inv:strong}.
As mentioned earlier, (EP) always has solutions, i.e.,
$\zer(A+B)\neq\varnothing$.
Apply
\cite[Theorem 26.14(i)\&(ii)]{BC2017} with $\beta=1/2$.

\ref{i:findg3}: Observe that
\begin{equation}\label{e:proj:dp}
P_{\Delta^{\perp}}=\bigg(\frac{\Id}{2}+T\bigg)(\Id-\bR).
\end{equation}
Indeed,
since that
$\big(\Id/2+T\big)^{-1}=\Id-\bR+2P_{\Delta}$ by Fact~\ref{f:shift}\ref{i:r3},
that $\Id/2+T$ is linear, and $\Delta\subseteq\ker T$, we have
\begin{align}
\Id &=\bigg(\frac{\Id}{2}+T\bigg)\bigg(\frac{\Id}{2}+T\bigg)^{-1}
=\bigg(\frac{\Id}{2}+T\bigg)(\Id-\bR)+2\bigg(\frac{\Id}{2}+T\bigg)P_{\Delta}\\
&=\bigg(\frac{\Id}{2}+T\bigg)(\Id-\bR)+P_{\Delta},
\end{align}
from which
$$P_{\Delta^{\perp}}=\Id-P_{\Delta}=\bigg(\frac{\Id}{2}+T\bigg)(\Id-\bR).$$
This, together with \ref{i:findg2}, Corollary~\ref{c:pd:sol}\ref{i:bx} and Theorem~\ref{t:many}, gives that when $n\rightarrow\infty,$
\begin{align}
P_{\Delta^{\perp}}(\bx_{n}) &=\bigg(\frac{\Id}{2}+T\bigg)(\bx_{n}-\bR\bx_{n})
\rightarrow \bigg(\frac{\Id}{2}+T\bigg)(-\by) =\bx=P_{\Delta^{\perp}}(\tilde{\bx}).
\end{align}
\end{proof}

When $\lambda_{n}\equiv 1$, an immediate consequence of Theorem~\ref{t:solve} comes as follows.
\begin{corollary} Let $\gamma\in\, ]0, 1[$, let $\bx_{0}\in \bX$ and set
\begin{equation}\label{e:recur}
(\forall n\in\NN)\ \bx_{n+1}=\prox_{\gamma\bg}\big((1-\gamma)\bx_{n}+\gamma\bR \bx_{n}\big).
\end{equation}
Then $\bx_{n}\rightharpoonup \tilde{\bx}$, a generalized cycle of $\bg$;
$\bR x_{n}-\bx_{n}\rightarrow \by$, the unique generalized gap vector of $\bg$;
$P_{\Delta^{\perp}}\bx_{n}\rightarrow\bx$, the unique solution to $(\tilde{P})$.
\end{corollary}

\begin{remark} In this regard, see \cite{abjw2020} for finding cycles and gap vectors of compositions of projections, and
also \cite[Section 3.3.3]{BCR} for an abstract framework.
\end{remark}


\section*{Acknowledgments}
The authors thank two anonymous referees for constructive and insightful comments.
HHB and XW were supported by the Natural Sciences and
Engineering Research Council of Canada.


\begin{thebibliography}{999}
\sepp

\bibitem{ABRW} 
S.\ Alwadani, H.H.\ Bauschke, J.P.\ Revalski, and X.\ Wang,
Resolvents and Yosida approximations of displacement mappings of
isometries, \emph{Set-Valued and Variational Analysis}, in press, 2021.
\url{https://arxiv.org/abs/2006.04860}.


\bibitem{abjw2020} S.\ Alwadani, H.B.\ Bauschke, J.P.\ Revalski, and X. Wang,
The difference vectors for convex sets and a resolution of the geometry conjecture, December 2020,
submitted, \url{https://arxiv.org/pdf/2012.04784.pdf}


\bibitem{AtTh} 
H.\ Attouch and M.\ Th\'era,
A general duality principle for the sum of two operators,
\emph{Journal of Convex Analysis}~3 (1996), 1--24.


\bibitem{bbr78}J.B.\ Baillon, R.E.\ Bruck, and S. Reich,
On the asymptotic behavior of nonexpansive mappings and semigroups in Banach spaces,
\emph{Houston Journal of Mathematics}~4 (1978), 1--9.

\bibitem{baillon}
J.-B.\ Baillon, P.L.\ Combettes, and R.\ Cominetti,
There is no variational characterization of the cycles in the method of periodic projections,
\emph{Journal of Functional Analysis}~262 (2012), 400--408.

\bibitem{baillon1}J.-B.\ Baillon, P.L.\ Combettes, and R.\ Cominetti,
Asymptotic behavior of compositions of under-relaxed nonexpansive operators,
\emph{Journal of Dynamics and Games}~1 (2014), 331--346.

\bibitem{bot}H.H.\ Bauschke, R.I.\ Bot, W.L.\ Hare, and W.\ Moursi,
Attouch-Th\'era duality revisited: paramonotonicity and operator splitting,
\emph{Journal of Approximation Theory}~164
(2012), 1065--1084.

\bibitem{BC2017}
H.H.\ Bauschke and P.L.\ Combettes,
\emph{Convex Analysis and Monotone Operator Theory in Hilbert Spaces},
second edition,
Springer, 2017.

 \bibitem{BCR}
 H.H.\ Bauschke, P.L.\ Combettes, and S.\ Reich,
 The asymptotic behavior of the composition of two resolvents.
 \emph{Nonlinear Analysis}~60 (2005), 283--301.

\bibitem{Victoria}
H.H.\ Bauschke, V.\ Mart\'in-M\'arquez, S.M.\ Moffat, and X.\ Wang,
Compositions and convex combinations of asymptotically regular
firmly nonexpansive mappings are also asymptotically regular,
\emph{Fixed Point Theory and Applications}~2012:53.

\bibitem{reich}
H.H.\ Bauschke, E.\ Matou\~{s}kov\'{a}, and S.\ Reich,
Projection and proximal point methods: convergence results and counterexamples,
\emph{Nonlinear Analysis}~56 (2004), 715--738.

\bibitem{bwy14}
H.H.\ Bauschke, X.\ Wang, and L.\ Yao,
Rectangularity and paramonotonicity of maximally monotone operators,
\emph{Optimization}~63 (2014), 487--504.

\bibitem{moffat13} H.H.\ Bauschke, S.M.\ Moffat, and X.\ Wang,
Near equality, near convexity, sums of maximally monotone operators, and averages of firmly nonexpansive mappings,
\emph{Mathematical Programming}~139 (2013), Ser. B, 55--70.

\bibitem{bruck} R.E.\ Bruck and S.\ Reich,
Nonexpansive projections and resolvents of accretive operators in Banach spaces,
\emph{Houston Journal of Mathematics}~3 (1977), 459--470.

\bibitem{comb20} P.L.\ Combettes and J.-C.\ Pesquet,
Fixed point strategies in data science, \emph{IEEE Transactions on Signal Processing},
doi: 10.1109/TSP.2021.3069677, 2021. \url{https://arxiv.org/pdf/2008.02260.pdf}

\bibitem{com-pes}
P.L.\ Combettes and J.-C.\ Pesquet,
Deep neural network structures solving variational inequalities,
\emph{Set-Valued and Variational Analysis}~28 (2020), 491--518.



%
%

%
%
%
%
%
%
%
 \bibitem{GR}
 K.\ Goebel and S.\ Reich,
 \emph{Uniform Convexity, Hyperbolic Geometry, and Nonexpansive
 Mappings}, Marcel Dekker, New York, 1984.
%
%
%
%
%

\bibitem{meyer2000} C.D.\ Meyer, \emph{Matrix Analysis and Applied Linear Algebra},
SIAM, 2000.

\bibitem{reich1}
S. Reich, On infinite products of resolvents,
\emph{Atti della Accademia Nazionale dei Lincei. Rendiconti. Classe di Scienze Fisiche, Matematiche e Naturali. Serie VIII}~63 (1977), 338--340.

\bibitem{reich3} S.\ Reich,
On the asymptotic behavior of nonlinear semigroups and the range of accretive operators,
\emph{Journal of Mathematical Analysis and Applications}~79 (1981), 113--126.

\bibitem{reich2} S.\ Reich and A.J.\ Zaslavski, Infinite products of resolvents of accretive operators,
\emph{Topological Methods in Nonlinear Analysis}~15 (2000), 153--168.

 \bibitem{Rock70}
 R.T.\ Rockafellar,
 \emph{Convex Analysis},
 Princeton University Press, Princeton, 1970.

%
%




\bibitem{Rock98}
R.T.\ Rockafellar and R.J-B\ Wets,
\emph{Variational Analysis},
Springer-Verlag, 
corrected 3rd printing, 2009.


\bibitem{Simons2}
S.\ Simons,
\emph{From Hahn-Banach to Monotonicity},
Springer-Verlag,
2008.

%
%
%
%

\bibitem{suzuki} T.\ Suzuki,
Some notes on Bauschke's condition, \emph{Nonlinear Analysis}~67 (2007), 2224--2231.

\bibitem{wangb} X.\ Wang and H.H.\ Bauschke,
Compositions and averages of two resolvents: relative geometry of fixed points sets and a partial
answer to a question by C. Byrne, \emph{Nonlinear Analysis}~74 (2011), 4550--4572.


 \bibitem{Zalinescu}{C.\ Z\u{a}linescu},
 \emph{Convex Analysis in General Vector Spaces},
 World Scientific Publishing, 2002.


\end{thebibliography}
\end{document}